\documentclass[11pt,reqno]{amsart}
\usepackage{amssymb,amscd}
\newtheorem{thm}{Theorem}[section]
\newtheorem{prop}[thm]{Proposition}
\newtheorem{lemma}[thm]{Lemma}
\newtheorem{cor}[thm]{Corollary}
\newtheorem{remark}[thm]{Remark}

\def\U{\mathrm{U}}
\def\bC{\mathbb{C}}
\def\cM{\mathcal{M}}
\def\bN{\mathbb{N}}
\def\eps{\varepsilon}
\def\proj{\mathrm{proj}}

\def\1{\mathbf{1}}
\def\cA{\mathcal{A}}

\def\ffi{\varphi}
\def\bR{\mathbb{R}}

\def\<{\langle}
\def\>{\rangle}

\def\Tr{\mathrm{Tr}}
\def\bZ{\mathbb{Z}}
\def\bT{\mathbb{T}}

\begin{document}

\title[Microstate free entropy of projections]
{Notes on microstate free entropy of projections}
\author[F. Hiai]{Fumio Hiai$\,^{1,2}$}
\address{Graduate School of Information Sciences,
Tohoku University, Aoba-ku, Sendai 980-8579, Japan}
\author[Y. Ueda]{Yoshimichi Ueda$\,^{1,3}$}
\address{Graduate School of Mathematics,
Kyushu University, Fukuoka 810-8560, Japan}

\thanks{$^1\,$Supported in part by Japan Society for the Promotion of Science,
Japan-Hungary Joint Project.}
\thanks{$^2\,$Supported in part by Grant-in-Aid for Scientific Research
(B)17340043.}
\thanks{$^3\,$Supported in part by Grant-in-Aid for Young Scientists
(B)17740096.}

\maketitle

\begin{abstract}
We study the microstate free entropy $\chi_{\proj}(p_1,\dots,p_n)$ of projections, and
establish its basic properties similar to the self-adjoint variable case.
Our main contribution is to characterize the pair-block freeness of projections by the
additivity of $\chi_{\proj}$ (Theorem \ref{T-4.1}), in the proof of which
a transportation cost inequality plays an important role. We also briefly discuss the
free pressure in relation to $\chi_{\proj}$.   
\end{abstract}

\section*{Introduction}

The theory of free entropy, initiated and mostly developed by D.~Voiculescu in his
series of papers \cite{V1}--\cite{V6}, has become one of the most essential disciplines
of free probability theory. For self-adjoint non-commutative random variables, say
$X_1,\dots,X_n$, the microstate free entropy  $\chi(X_1,\dots,X_n)$ introduced in
\cite{V2} is defined as a certain asymptotic growth rate (as the matrix size $N$ goes
to $\infty$) of the Euclidean volume of the set of $N\times N$  self-adjoint matrices
$(A_1,\dots,A_n)$ approximating $(X_1,\dots,X_n)$ in moments. It is this microstate
theory that settled some long-standing open questions in von Neumann algebras (see the
survey \cite{V-Survey}). On the other hand, the non-microstate free entropy
$\chi^*(X_1,\dots,X_n)$ was also introduced in \cite{V5} based on the non-commutative
Hilbert transform and the notion of conjugate variables, without the use
of microstates or so-called matrix integrals which are rather hard to handle.
Although it is believed that both approaches should be unified and give
the same quantity, only the inequality $\chi\le\chi^*$ is known to hold true due to
Biane, Capitaine and Guionnet \cite{BCG} based on an idea of large deviation principle
for several random matrices. In his work \cite{V6} Voiculescu developed
another kind of non-microstate approach to the free entropy, the so-called free
liberation theory, and introduced the mutual free information $i^*(X_1,\dots,X_n)$
based on it. He suggested there the need to apply the microstate approach
to projection random variables because the usual microstate free entropy $\chi$ becomes
always zero for projections while $i^*$ does not. Following the suggestion, we
here study the microstate free entropy $\chi_\proj(p_1,\dots,p_n)$ of
projections $p_1,\dots,p_n$ in the same lines as in \cite{V2} and \cite{V4}
to provide the basis for future research.

The large deviation principle for random matrices as mentioned above
started with the paper of Ben Arous and Guionnet \cite{BG} and has been almost
completed in the single random matrix case (corresponding to the study of $\chi(X)$ for
single random variable $X$), see the survey \cite{Gu-Survey}. We note that such large
deviation principle played quite an important role not only for the foundation of free
entropy theory but also for getting free analogs of several probability theoretic
inequalities (see \cite{HU1} and the references therein). Recently, one
more large deviation was shown in \cite{HP1} for an independent pair of random
projection matrices, including the explicit formula of the free entropy
$\chi_\proj(p,q)$ of a projection pair $(p,q)$. This is one of a few large deviation
results (indeed the first full large deviation result) in the setting of several random
matrices, though the method of the proof is based on the single variable
case. Moreover, in \cite{HU2} we applied it to get a kind of logarithmic Sobolev
inequality $\chi_\proj(p,q)\le\ffi^*(p:q)$ between the free entropy $\chi_\proj(p,q)$
and the mutual free Fisher information $\ffi^*(p,q)$ (see \cite{V6}) for a projection
pair. The large deviation result in \cite{HP1} also plays a crucial role
in our study of $\chi_{\proj}$ here.

The paper is organized as follows. After stating the definition and basic properties
of $\chi_\proj(p_1,\dots,p_n)$ in \S1, we recall in \S2 the formula in the case of two
variables. In \S3 we introduce a certain functional calculus for a projection pair
$(p,q)$ and provide a technical tool of separate change of variable formula. This
tool is essential in \S4 to prove the additivity theorem characterizing the
pair-block freeness of projections by the additivity of their free entropy.
\S5 treats a free analog of transportation cost inequalities for tracial distributions
of projections. Its simplest case is needed in the proof of the above additivity
theorem while of interest by itself. Finally, along the same lines as in \cite{Hi}, we
introduce in \S6 the notion of free pressure and compare its Legendre transform with
$\chi_\proj(p_1,\dots,p_n)$, thus giving a variational expression of free entropy.

\section{Definition}
\setcounter{equation}{0}

For $N\in\bN$ let $\U(N)$ be the unitary group of order $N$. For $k\in\{0,1,\dots,N\}$
let $G(N,k)$ denote the set of all $N\times N$ orthogonal projection matrices of rank
$k$, that is, $G(N,k)$ is identified with the Grassmannian manifold consisting of
$k$-dimensional subspaces in $\bC^N$. With the diagonal matrix $P_N(k)$ of the first
$k$ diagonals $1$ and the others $0$, each $P\in G(N,k)$ is diagonalized as
\begin{equation}\label{F-1.1}
P=UP_N(k)U^*,
\end{equation}
where $U\in\U(N)$ is determined up to the right multiplication of elements in
$\U(k)\oplus\U(N-k)$. Hence $G(N,k)$ is identified with the homogeneous space
$\U(N)/(\U(k)\oplus\U(N-k))$, and we have a unique probability measure
$\gamma_{G(N,k)}$ on $G(N,k)$ invariant under the unitary conjugation $P\mapsto UPU^*$
for $U \in \U(N)$. Via the description as homogeneous space, this
corresponds to the measure on $\U(N)/(\U(k)\oplus\U(N-k))$ invarinat under the left
multiplication of elements in $\U(N)$ or induced from the Haar probability measure
$\gamma_{\U(N)}$ on $\U(N)$. Let $\xi_{N,k}:\U(N)\to G(N,k)$ be the (surjective
continuous) map defined by \eqref{F-1.1}, i.e., $\xi_{N,k}(U):=UP_N(k)U^*$. Then the
measure $\gamma_{G(N,k)}$ is more explicitly written as
\begin{equation}\label{F-1.2}
\gamma_{G(N,k)}=\gamma_{\U(N)}\circ\xi_{N,k}^{-1}.
\end{equation}

Throughout the paper $(\cM,\tau)$ is a tracial $W^*$-probability space. Let
$(p_1,\dots,p_n)$ be an $n$-tuple of projections in $(\cM,\tau)$. Following
Voiculescu's proposal in \cite[14.2]{V6} we define the {\it free entropy}
$\chi_\proj(p_1,\dots,p_n)$ of $(p_1,\dots,p_n)$ as follows. Choose
$k_i(N)\in\{0,1,\dots,N\}$ for each $N\in\bN$ and $1\le i\le n$ in such a way that
$k_i(N)/N\to\tau(p_i)$ as $N\to\infty$ for $1\le i\le n$. For each $m\in\bN$ and
$\eps>0$ set
\begin{align}\label{F-1.3}
&\Gamma_\proj(p_1,\dots,p_n;k_1(N),\dots,k_n(N);N,m,\eps) \nonumber\\
&\quad:=\biggl\{(P_1,\dots,P_n)\in\prod_{i=1}^nG(N,k_i(N)):
\bigg|{1\over N}\Tr_N(P_{i_1}\cdots P_{i_r})-\tau(p_{i_1}\cdots p_{i_r})\bigg|<\eps
\nonumber\\
&\hskip6cm\mbox{for all $1\le i_1,\dots,i_r\le n$, $1\le r\le m$}\biggr\},
\end{align}
where $\Tr_N$ stands for the usual (non-normalized) trace on the $N\times N$ matrices. 
We then define
\begin{align}\label{F-1.4}
&\chi_\proj(p_1,\dots,p_n):=
\lim_{\substack{m\rightarrow\infty \\ \varepsilon \searrow 0}}
\limsup_{N\to\infty} \nonumber\\
&\qquad{1\over N^2}\log\Biggl(\bigotimes_{i=1}^n\gamma_{G(N,k_i(N))}\Biggr)
\bigl(\Gamma(p_1,\dots,p_n;k_1(N),\dots,k_n(N);N,m,\eps)\bigr).
\end{align}
For the justification of the definition of $\chi_{\proj}$, here arises a natural
question whether the quantity $\chi_{\proj}(p_1,\dots,p_n)$ depends on the particular
choice of $k_i(N)$ or not. The following is the answer to it.

\begin{prop}\label{P-1.1}
The above definition of $\chi_\proj(p_1,\dots,p_n)$ is independent of the choices of
$k_i(N)$ with $k_i(N)/N\to\alpha_i$ for $1\le i\le n$.
\end{prop}

\begin{proof}
For $1\le i\le n$ let $l_i(N)$, $N\in\bN$, be another sequence such that
$l_i(N)/N\to\alpha_i$ as $N\to\infty$. For each $N,m\in\bN$ and $\eps>0$, we write
$\Gamma(\vec k(N),m,\eps)$ ($\subset\prod_{i=1}^nG(N,k_i(N))$) for the set
\eqref{F-1.3} with respect to $\vec k(N):=(k_1(N),\dots,k_n(N))$, and also
$\Gamma(\vec l(N),m,\eps)$ ($\subset\prod_{i=1}^nG(N,l_i(N))$) for the same with
respect to $\vec l(N):=(l_1(N),\dots,l_n(N))$. Moreover, we set
$\xi_{\vec k(N)}(\vec U):=\bigl(\xi_{N,k_1(N)}(U_1),\dots,\xi_{N,k_n(N)}(U_n)\bigr)$
for $\vec U=(U_1,\dots,U_n)\in\U(N)^n$, and define the subset
$\widetilde\Gamma(\vec l(N),m,\eps):=\xi_{\vec l(N)}\circ
\xi_{\vec k(N)}^{-1}\bigl(\Gamma(\vec k(N),m,\eps)\bigr)$
of $\prod_{i=1}^nG(N,l_i(N))$. For every $N\in\bN$ and $U\in\U(N)$, since
$$
\xi_{N,l_i(N)}(U)-\xi_{N,k_i(N)}(U)=U\bigl(P_N(k_i(N))-P_N(l_i(N))\bigr)U^*,
$$
we get
$$
\big\|\xi_{N,l_i(N)}(U)-\xi_{N,k_i(N)}(U)\big\|_1={|l_i(N)-k_i(N)|\over N},
$$
where $\|\cdot\|_1$ denotes the trace-norm with respect to $N^{-1}\Tr_N$. For every
$m\in\bN$ and $\eps>0$, there exists an $N_0\in\bN$ such that
$N^{-1}|l_i(N)-k_i(N)|<\eps/m$ for all $N\ge N_0$ and $1\le i\le n$. Let us prove that
$\widetilde\Gamma(\vec l(N),m,\eps)\subset\Gamma(\vec l(N),m,2\eps)$ whenever
$N\ge N_0$. Assume that $N\ge N_0$ and
$\vec Q=(Q_1,\dots,Q_n)\in\widetilde\Gamma(\vec l(N),m,\eps)$; then
$\vec U=(U_1,\dots,U_n)\in\U(N)^n$ exists so that $\vec Q=\xi_{\vec l(N)}(\vec U)$ and
$\vec P=(P_1,\dots,P_n):=\xi_{\vec k(N)}(\vec U)\in\Gamma_{\vec k(N),m,\eps}$. Since
$$
\|Q_i-P_i\|_1=\big\|\xi_{N,l_i(N)}(U_i)-\xi_{N,k_i(N)}(U_i)\big\|_1
<{\eps\over m},\qquad1\le i\le n,
$$
we get for $1\le i_1,\dots,i_r\le n$ and $1\le r\le m$
$$
\bigg|{1\over N}\Tr_N(Q_{i_1}\cdots Q_{i_r})
-{1\over N}\Tr_N(P_{i_1}\cdots P_{i_r})\bigg|
\le\sum_{j=1}^r\|Q_{i_j}-P_{i_j}\|_1<\eps
$$
so that
$$
\bigg|{1\over N}\Tr_N(Q_{i_1}\cdots Q_{i_r})-\tau(p_{i_1}\cdots p_{i_r})\bigg|<2\eps,
$$
implying $\vec Q\in\Gamma(\vec l(N),m,2\eps)$. Setting
$\gamma_{\vec k(N)}:=\bigotimes_{i=1}^n\gamma_{G(N,k_i(N))}$, we now have thanks to
\eqref{F-1.2}
\begin{align*}
\gamma_{\vec l(N)}\bigl(\Gamma(\vec l(N),m,2\eps)\bigr)
&\ge\gamma_{\vec k(N)}\bigl(\widetilde\Gamma(\vec l(N),m,\eps)\bigr) \\
&=\bigl(\gamma_{\U(N)}\bigr)^{\otimes n}
\circ\xi_{\vec l(N)}^{-1}\circ\xi_{\vec l(N)}\circ\xi_{\vec k(N)}^{-1}
\bigl(\Gamma(\vec k(N),m,\eps)\bigr) \\
&\ge\bigl(\gamma_{\U(N)}\bigr)^{\otimes n}\circ\xi_{\vec k(N)}^{-1}
\bigl(\Gamma(\vec k(N),m,\eps)\bigr)
=\gamma_{\vec k(N)}\bigl(\Gamma(\vec k(N),m,\eps)\bigr)
\end{align*}
whenever $N\ge N_0$. This implies that
$$
\limsup_{N\to\infty}{1\over N^2}\log
\gamma_{\vec l(N)}\bigl(\Gamma(\vec l(N),m,2\eps)\bigr) \\
\ge\limsup_{N\to\infty}{1\over N^2}\log
\gamma_{\vec k(N)}\bigl(\Gamma(\vec k(N),m,\eps)\bigr),
$$
which says that the free entropy \eqref{F-1.4} given for $\vec k(N)$ is not greater
than that for $\vec l(N)$. By symmetry we observe that both free entropies must
coincide.
\end{proof}

The following are basic properties of $\chi_\proj$. We omit their proofs, all of which
are essentially same as in the case of self-adjoint variables in \cite{V2} or else obvious.

\begin{prop}\label{P-1.2}
Let $p_1,\dots,p_n$ be projections in $(\cM,\tau)$.
\begin{itemize}
\item[(i)] Negativity{\rm :} $\chi_\proj(p_1,\dots,p_n)\le0$.
\item[(ii)] Subadditivity{\rm :} for every $1\le j<n$,
$$
\chi_\proj(p_1,\dots,p_n)\le\chi_{\proj}(p_1,\dots,p_j)+\chi_\proj(p_{j+1},\dots,p_n).
$$
\item[(iii)] Upper semi-continuity{\rm :} if a sequence $(p_1^{(m)},\dots,p_n^{(m)})$ of
$n$-tuples of projections converges to $(p_1,\dots,p_n)$ in distribution, then
$$
\chi_\proj(p_1,\dots,p_n)\ge\limsup_{m\to\infty}
\chi_\proj(p_1^{(m)},\dots,p_n^{(m)}).
$$
\item[(iv)] $\chi_\proj(p_1,\dots,p_n)$ does not change when $p_i$ is replaced by
$p_i^\perp:=\1-p_i$ for each $i$.
\end{itemize}
\end{prop}

\begin{remark}\label{R-1.3} {\rm
We may adopt different ways to introduce the free entropy of an $n$-tuple
$(p_1,\dots,p_n)$ of projections in $(\cM,\tau)$. For instance, for each $N\in\bN$
consider two unitarily invariant probability measures $\gamma_{G(N)}^{(1)}$ and
$\gamma_{G(N)}^{(2)}$ on $G(N):=\bigsqcup_{k=0}^NG(N,k)$ determined by the weights on
$G(N,k)$, $0\le k\le N$, given as
$$
\gamma_{G(N)}^{(1)}(G(N,k))={1\over N+1},\quad
\gamma_{G(N)}^{(2)}(G(N,k))={1\over2^N}{N\choose k}.
$$
For each $m\in\bN$ and $\eps>0$ set
\begin{align*}
&\Gamma_\proj(p_1,\dots,p_n;N,m,\eps) \\
&\qquad:=\biggl\{(P_1,\dots,P_n)\in G(N)^n:
\bigg|{1\over N}\Tr_N(P_{i_1}\cdots P_{i_r})-\tau(p_{i_1}\cdots p_{i_r})\bigg|<\eps \\
&\hskip6cm\mbox{for all $1\le i_1,\dots,i_r\le n$, $1\le r\le m$}\biggr\},
\end{align*}
and define for $j=1,2$
$$
\chi_\proj^{(j)}(p_1,\dots,p_n)
:=\lim_{\substack{m\rightarrow\infty \\ \varepsilon \searrow 0}}
\limsup_{N\to\infty}{1\over N^2}
\log\Bigl(\gamma_{G(N)}^{(j)}\Bigr)^{\otimes n}
\bigl(\Gamma_\proj(p_1,\dots,p_n;N,m,\eps)\bigr).
$$
It is fairly easy to see (similarly to the proof of Proposition 1.1) that
both $\chi_\proj^{(j)}(p_1,\dots,p_n)$, $j=1,2$, coincide with
$\chi_\proj(p_1,\dots,p_n)$ given in \eqref{F-1.4}.
}\end{remark}

\section{Case of two projections}
\setcounter{equation}{0}

Let $(p,q)$ be a pair of projections in $(\cM,\tau)$ with $\alpha:=\tau(p)$ and
$\beta:=\tau(q)$. Set
$$
E_{11}:=p\wedge q,\quad E_{10}:=p\wedge q^\perp,\quad
E_{01}:=p^\perp\wedge q,\quad E_{00}:=p^\perp\wedge q^\perp,
$$
$$
E:=\1-(E_{00}+E_{01}+E_{10}+E_{11}).
$$
Then $E$ and $E_{ij}$ are in the center of $\{p,q\}''$ and
$(E\{p,q\}''E,\tau|_{E\{p,q\}''E})$ is isomorphic to $L^\infty((0,1),\nu;M_2(\bC))$,
where $\nu$ is the measure on $(0,1)$ determined by 
$$
\tau(A) = \frac{1}{2} \int_{(0,1)} \mathrm{Tr}_2(A(x))\,d\nu(x),
\quad A \in L^{\infty}((0,1),\nu;M_2(\mathbb{C})) \cong E\{p,q\}'' E
$$
(hence $\nu((0,1))=\tau(E)$). Under this isomorphism, $EpE$ and $EqE$
are represented as
$$
(EpE)(x)=\bmatrix1&0\\0&0\endbmatrix\ \ \mbox{and}
\ \ (EqE)(x)=\bmatrix x&\sqrt{x(1-x)}\\\sqrt{x(1-x)}&1-x\endbmatrix
\quad\mbox{for $x\in(0,1)$}. 
$$
In this way, the mixed moments of $(p,q)$ with respect to $\tau$ are determined by
$\nu$ and $\{\tau(E_{ij})\}_{i,j=0}^1$. Although $\nu$ is not necessarily a
probability measure, we define the free entropy $\Sigma(\nu)$ by
$$
\Sigma(\nu):=\iint_{(0,1)^2}\log|x-y|\,d\nu(x)\,d\nu(y)
$$
in the same way as in \cite{V1}. Furthermore, we set
\begin{equation}\label{F-2.1}
\rho:=\min\{\alpha,\beta,1-\alpha,1-\beta\},
\end{equation}
\begin{equation}\label{F-2.2}
C:=\rho^2B\biggl({|\alpha-\beta|\over\rho},{|\alpha+\beta-1|\over\rho}\biggr)
\end{equation}
(meant zero if $\rho=0$), where
\begin{align*}
B(s,t)&:={(1+s)^2\over2}\log(1+s)-{s^2\over2}\log s
+{(1+t)^2\over2}\log(1+t)-{t^2\over2}\log t \\
&\qquad-{(2+s+t)^2\over2}\log(2+s+t)+{(1+s+t)^2\over2}\log(1+s+t)
\end{align*}
for $s,t\geq0$. With these definitions, the following formula of $\chi_\proj(p,q)$ was
obtained in \cite{HP1} as a consequence of the large deviation principle for an
independent pair of random projection matrices.

\begin{prop}\label{P-2.1}{\rm (\cite[Theorem 3.2, Proposition 3.3]{HP1})}\quad
If $\tau(E_{00})\tau(E_{11})=\tau(E_{01})\tau(E_{10})=0$, then
\begin{align*}
\chi_\proj(p,q)&={1\over4}\Sigma(\nu)
+{|\alpha-\beta|\over2}\int_{(0,1)}\log x\,d\nu(x) \\
&\qquad\quad+{|\alpha+\beta-1|\over2}\int_{(0,1)}\log(1-x)\,d\nu(x)-C,
\end{align*}
and otherwise $\chi_\proj(p,q)=-\infty$. Moreover, $\chi_\proj(p,q)=0$ if and only if
$p$ and $q$ are free.
\end{prop}

Note that the condition $\tau(E_{00})\tau(E_{11})=\tau(E_{01})\tau(E_{10})=0$ is 
equivalent to
\begin{equation}\label{F-2.3}
\begin{cases}
\tau(E_{11})=\max\{\alpha+\beta-1,0\}, \\
\tau(E_{00})=\max\{1-\alpha-\beta,0\}, \\
\tau(E_{10})=\max\{\alpha-\beta,0\}, \\
\tau(E_{01})=\max\{\beta-\alpha,0\};
\end{cases}
\end{equation}
in this case, $\tau(E_{01})+\tau(E_{10})=|\alpha-\beta|$,
$\tau(E_{00})+\tau(E_{11})=|\alpha+\beta-1|$ and $\tau(E)=2\rho$.

In the case where $\tau_\proj(p,q)=0$ (equivalently, $p$ and $q$ are free), the
measure $\nu$ was computed in \cite{VDN} as
\begin{equation}\label{F-2.4}
{\sqrt{(x-\xi)(\eta-x)}\over2\pi x(1-x)}\1_{(\xi,\eta)}(x)\,dx
\end{equation}
with $\xi,\eta:=\alpha+\beta-2\alpha\beta\pm\sqrt{4\alpha\beta(1-\alpha)(1-\beta)}$.
It is also worthwhile to note (see \cite{HP1}) that $\limsup$ in definition
\eqref{F-1.4} can be replaced by $\lim$ in the case of two projections.

In \S4 the equivalence between the additivity of $\chi_\proj$ and the freeness of
projection variables will be generalized to the case of more than two projections. To
do this, we need a kind of separate change of variable formula for $\chi_\proj$
established in the next section. 

\section{Separate change of variable formula}
\setcounter{equation}{0}

Let $N\in\bN$ and $k,l\in\{0,1,\dots,N\}$. Assume that $0<k\le l$ and $k+l\le N$.
Consider a pair $(P,Q)$ of $N\times N$ projection matrices with ${\rm rank}(P)=k$ and
${\rm rank}(Q)=l$, which is distributed under the measure
$\gamma_{G(N,k)}\otimes\gamma_{G(N,l)}$ on $G(N,k)\times G(N,l)$. Thanks to
the assumptions on $k,l$, for any pair $(P,Q)\in G(N,k)\times G(N,l)$ the so-called
sine-cosine decomposition of two projections gives the following representation:
\begin{align}
P&=U\left(\bmatrix I&0\\0&0\endbmatrix\oplus0\oplus0\right)U^*, \label{F-3.1}\\
Q&=U\left(\bmatrix X&\sqrt{X(I-X)}\\\sqrt{X(I-X)}&I-X\endbmatrix
\oplus I\oplus0\right)U^* \label{F-3.2}
\end{align}
in $\bC^N=(\bC^k\otimes\bC^2)\oplus\bC^{l-k}\oplus\bC^{N-k-l}$, 
where $U$ is an $N\times N$ unitary matrix and $X$ is a $k\times k$ diagonal matrix
with the diagonal entries $0\leq x_1\leq x_2\leq \dots\leq x_k\leq 1$. When
$x_1,\dots,x_k$ are in $(0,1)$ and mutually distinct, it is easy to see that $U$ is
uniquely determined up to the right multiplication of unitary matrices of the form
$$
\bmatrix T&0\\0&T\endbmatrix\oplus V_1\oplus V_2,
\quad T\in\bT^k,\ V_1\in\U(l-k),\ V_2\in\U(N-k-l).
$$
We denote by $V(N,k,l)$ the subgroup of $\U(N)$ consisting of all unitary matrices of
the above form so that $\U(N)/V(N,k,l)$ becomes a homogeneous space. Also, let
$[0,1]_{\leq}^k$ and $(0,1)_<^k$ denote the sets of $(x_1,\dots,x_k)$ satisfying
$0\leq x_1\leq\dots\leq x_N\leq 1$ and $0<x_1<\dots<x_k<1$, respectively. We then
consider the continuous map $\Xi_{N,k,l}:
\mathrm{U}(N)/V(k,\ell) \times [0,1]_{\leq}^k \rightarrow G(N,k)\times G(N,l)$
defined by \eqref{F-3.1} and \eqref{F-3.2}, that is,
\begin{align*}
&\Xi_{N,k,l}([U],X) \\
&\quad := 
\left(
U\left(\begin{bmatrix} I & 0 \\ 0 & 0 \end{bmatrix} \oplus0\oplus0\right)U^*, 
U\left(\begin{bmatrix} X & \sqrt{X(I-X)} \\ \sqrt{X(I-X)} & I-X \end{bmatrix}
\oplus I\oplus0\right)U\right), 
\end{align*} 
where $X$ is regarded as a diagonal matrix. The set
$$
\left(G(N,k)\times G(N,l)\right)_0:=
\Xi_{N,k,l}\left(\mathrm{U}(N)/V(N,k,l)\times (0,1)^k_<\right)
$$
is open and co-negligible with respect to $\gamma_{G(N,k)}\otimes\gamma_{G(N,l)}$
in $G(N,k)\times G(N,l)$ thanks to \cite[Theorem 2.2]{Col} (or \cite[Lemma 1.1]{HP1})
and moreover $\Xi_{N,k,l}$ gives a smooth diffeomorphism between
$(\U(N)/V(N,k,l))\times(0,1)_<^k$ and $\left(G(N,k)\times G(N,l)\right)_0$. The next
lemma will be needed later.

\begin{lemma}\label{L-3.1}
The measure $(\gamma_{G(N,k)}\otimes\gamma_{G(N,l)})\circ\Xi_{N,k,l}$ coincides
with 
$$
\gamma_{N,k,l}\otimes\Biggl({1\over Z_{N,k,l}}\prod_{i=1}^kx_i^{l-k}
(1-x_i)^{N-k-l}\prod_{1\le i<j\le k}(x_i-x_j)^2\prod_{i=1}^kdx_i\Biggr),
$$
where $\gamma_{N,k,l}$ is the {\rm (}unique{\rm )} probability measure on
$\U(N)/V(N,k,l)$ induced by the Haar probability measure on $\U(N)$ and $Z_{N,k,l}$ is
a normalization constant.
\end{lemma}

\begin{proof}
Let $\lambda$ be the measure on $\mathrm{U}(N)/V(N,k,l)\times(0,1)^k_<$ transformed
from the restriction of $\gamma_{G(N,k)}\otimes\gamma_{G(N,l)}$ to
$\left(G(N,k)\times G(N,l)\right)_0$ by the inverse of $\Xi_{N,k,l}$, and $\mu$
be its image measure by the projection map $([U],X) \mapsto X$. The
disintegration theorem (see e.g.~\cite[Chapter IV, \S6.5]{Ma}) ensures that there is a
$\mu$-a.e.~unique Borel map $\lambda_{(\cdot)}$ from $(0,1)^k_<$ to the probability
measures on $\U(N)/V(N,k,l)$ such that $\lambda = \int_{(0,1)^k_<} \lambda_X\,d\mu(X)$.
Note that $([U],X) \mapsto X$ splits into $\Xi_{N,k,l}$, $(P,Q) \mapsto PQP$ and
the map sending $PQP$ to the eigenvalues in increasing order. Hence $\mu$ coincides
with the eigenvalue distribution of $PQP$ arranged in increasing order, which is
known to be equal to the second component given in the lemma by
\cite[Theorem 2.2]{Col}. Therefore, it suffices to show that $\lambda_X$ coincides
with $\gamma_{N,k,l}$ for $\mu$-a.e.~$X \in (0,1)_<^k$. For each $V\in\U(N)$, the
unitary conjugation $\mathrm{Ad}\,V\times\mathrm{Ad}\,V : (P,Q)\mapsto(VPV^*,VQV^*)$
on $G(N,k)\times G(N,l)$ and the left-translation $L_V : [U]\mapsto V[U] := [VU]$ on
$\mathrm{U}(N)/V(N,k,l)$ satisfy the relation $\Xi_{N,k,l}\circ
(L_V\times\mathrm{id})= (\mathrm{Ad}\,V\times\mathrm{Ad}\,V)\circ\Xi_{N,k,l}$;
hence, in particular, $\left(G(N,k)\times G(N,l)\right)_0$ is invariant under the
action $\mathrm{Ad}\,V\times\mathrm{Ad}\,V$ for every $V \in \mathrm{U}(N)$. Then,
for any bounded Borel function $f$ on $\mathrm{U}(N)/V(N,k,l)\times(0,1)^k)$, one can
easily verify that
\begin{align*}
&\int_{(0,1)_<^k}\left(\int_{\mathrm{U}(N)/V(N,k,l)} f([U],X)
\,d(\lambda_X\circ L_V)([U])\right)\,d\mu(X) \\
&\qquad\quad=\int_{\mathrm{U}(N)/V(N,k,l)\times(0,1)_<^k} f([U],X)\,d\lambda([U],X),
\end{align*}
which means a new disintegration
$\lambda = \int_{(0,1)_<^k} \lambda_X\circ L_V\,d\mu(X)$. The uniqueness of the
disintegration says that for $\mu$-a.e.~$X \in (0,1)_<^k$ one has
$\lambda_X = \lambda_X\circ L_V$ for all $V \in \mathrm{U}(N)$. Since $\gamma_{N,k,l}$
is a unique probability measure on $\mathrm{U}(N)/V(N,k,l)$ invariant under the
left-translation action of $\mathrm{U}(N)$, it follows that $\lambda_X=\gamma_{N,k,l}$ for $\mu$-a.e.~$X \in (0,1)^k_<$ so that 
$$
\lambda=\int_{(0,1)_<^k}\gamma_{N,k,l}\,d\mu(X)=\gamma_{N,k,l}\otimes\mu,
$$
as required.
\end{proof}

For a pair $(p,q)$ of projections in $(\cM,\tau)$ we introduce a sort of functional
calculus via the representation explained in \S2 in the following way. Let $\psi$ be
a continuous increasing function $\psi$ from $(0,1)$ into itself. With the notations
in \S2 we define a projection $q(\psi;p)$ in $\{p,q\}''$ by
$$
q(\psi;p):=Eq(\psi;p)E+E_{00}+E_{01}+E_{10}+E_{11},
$$
$$
(Eq(\psi;p)E)(x):=\bmatrix\psi(x)&\sqrt{\psi(x)(1-\psi(x))}\\
\sqrt{\psi(x)(1-\psi(x))}&1-\psi(x)\endbmatrix\quad\mbox{for $x\in(0,1)$}.
$$
It is obvious that $\tau(q(\psi;p))=\tau(q)$. (The definition itself is possible for
general Borel function from $(0,1)$ into $[0,1]$ but the above case is enough for our
purpose.)  The aim of this section is to prove the following change of variable
formula for free entropy of projections.

\begin{thm}\label{T-3.2}
Let $p_1,q_1,\dots,p_n,q_n,r_{1},\dots,r_{n'}$ be projections in $(\cM,\tau)$ and
assume that $\chi_{\rm proj}(p_i,q_i)>-\infty$ for $1\le i\le n$. Let
$\psi_1,\dots,\psi_n$ be continuous increasing functions from $(0,1)$ into itself,
and $q_i(\psi_i;p_i)$ be the projection defined from $p_i$, $q_i$ and $\psi_i$ as
above for $1\le i\le n$. Then we have
\begin{align*}
&\chi_{\rm proj}(p_1,q_1(\psi_1;p_1),\dots,p_n,q_n(\psi_n;p_n), r_1,\dots,r_{n'}) \\
&\quad\ge\chi_{\rm proj}(p_1,q_1,\dots,p_n,q_n,r_{1},\dots,r_{n'})
+\sum_{i=1}^n\bigl\{\chi_{\rm proj}(p_i,q_i(\psi_i;p_i))
-\chi_{\rm proj}(p_i,q_i)\bigr\}.
\end{align*}
Moreover, if $\psi_1,\dots,\psi_n$ are strictly increasing, then equality holds true
in the above inequality.
\end{thm}

The proof goes on the essentially same lines as in \cite{V4} and it is divided into
two steps; one is to analyze the case when $\psi_1,\dots,\psi_n$ are all extended
to $C^1$-diffeomorphisms from $[0,1]$ onto itself and the other is to approximate,
in two stages, the given $\psi_1,\dots,\psi_n$ by
$C^{\infty}$-diffeomorphisms from $[0,1]$ onto itself in such a way that the
corresponding free entropies converge to those in question. As  the first step let us
prove the following special case of the theorem.

\begin{lemma}\label{L-3.3}
Let $p_1,q_1,\dots,p_n,q_n,r_{1},\dots,r_{n'}$ be as in Theorem \ref{T-3.2}. If
$\psi_1,\dots,\psi_n$ are $C^1$-diffeomorphisms from $[0,1]$ onto itself with
$\psi_i(0) = 0$, $\psi_i(1) = 1$ and moreover $\psi_i'(x)>0$ for all $x\in[0,1]$,
then the equality of Theorem \ref{T-3.2} holds true. 
\end{lemma}

Obviously, it suffices to show when $n=1$; hence we assume $n=1$ and write $p=p_1$,
$q=q_1$ and $\psi=\psi_1$ for brevity. Let $\nu$ and $\{E_{ij}\}_{i,j=0}^1$
be as in \S2 for $(p,q)$. By Propositions \ref{P-1.2}\,(iv) and \ref{P-2.1} we may
assume that $\tau(p) \leq \tau(q) \leq 1/2$ so that $E_{11} = E_{10} = 0$
by \eqref{F-2.3}. We may further assume that $p$ is non-zero; otherwise there is
nothing to do. With the polar decomposition $(1-p)qp = v_{p,q}\sqrt{pqp(p-pqp)}$,
we thus represent $p$, $q$ and $q(\psi;p)$ as follows:
\begin{align*} 
p &= v_{p,q}^* v_{p,q}, \\
q &= pqp + v_{p,q}\sqrt{pqp(p-pqp)}
+ \sqrt{pqp(p-pqp)}v_{p,q}^* + v_{p,q}(p-pqp)v_{p,q}^* \\
&\quad+ \Bigl(q-pqp - (1-p)qp - pq(1-p) - v_{p,q}(p-pqp)v_{p,q}^*\Bigr), \\
q(\psi;p) &= \psi(pqp) + v_{p,q}\sqrt{\psi(pqp)(p-\psi(pqp))} \\
&\quad + \sqrt{\psi(pqp)(p-\psi(pqp))}v_{p,q}^*
+ v_{p,q}(p-\psi(pqp))v_{p,q}^* \\
&\quad + \Bigl(q-pqp - (1-p)qp - pq(1-p) - v_{p,q}(p-pqp)v_{p,q}^*\Bigr), 
\end{align*} 
where $\psi(pqp)$ means the functional calculus of $pqp$. Choose two sequences
$k(N)$, $l(N)$ for $N\ge2$ in such a way that $0<k(N) \leq l(N) \leq N/2$ and
$k(N)/N \rightarrow \tau(p)$, $l(N)/N \rightarrow \tau(q)$ as $N \rightarrow \infty$. 
As explained at the beginning of this section, for each
$(P,Q) \in (G(N,k(N))\times G(N,l(N)))_0$ there is a unitary $U \in \mathrm{U}(N)$,
unique up to $V(N,k(N),l(N))$, for which we have \eqref{F-3.1} and \eqref{F-3.2}.
Then we can define the map $\Phi_{N,\psi}$ on $(G(N,k(N))\times G(N,l(N)))_0$ by
sending $(P,Q)$ to $(P,Q(\psi;P))$ with  
$$ 
Q(\psi;P) := U\left(\begin{bmatrix}\psi(X) & \sqrt{\psi(X)(I-\psi(X))} \\
\sqrt{\psi(X)(I-\psi(X))} & I-\psi(X) \end{bmatrix} \oplus I\oplus0\right)U^*. 
$$
With the polar decomposition $(I-P)QP = V_{P,Q}\sqrt{PQP(I-PQP)}$ we have the following
expressions:
\begin{align*} 
Q &= PQP + V_{P,Q}\sqrt{PQP(P-PQP)} \\
&\quad +\sqrt{PQP(P-PQP)}V_{P,Q}^* + V_{P,Q}(P-PQP)V_{P,Q}^* \\
&\quad+ \Bigl(Q-PQP-(I-P)QP-PQ(I-P)-V_{P,Q}(P-PQP)V_{P,Q}^*\Bigr), \\
Q(\psi;P) &= \psi(PQP) + V_{P,Q}\sqrt{\psi(PQP)(P-\psi(PQP))} \\
&\quad +\sqrt{\psi(PQP)(P-\psi(PQP))}V_{P,Q}^* + V_{P,Q}(P-\psi(PQP))V_{P,Q}^* \\
&\quad + \Bigl(Q-PQP-(I-P)QP-PQ(I-P)-V_{P,Q}(P-PQP)V_{P,Q}^*\Bigr).  
\end{align*}   
Upon these expressions, what we now need is to approximate $v_{p,q}$ 
and $V_{P,Q}$ by polynomials of $p,q$ and $P,Q$, respectively, as stated in the
next lemma very similarly to \cite[Lemma 2.6]{HP1} (or \cite[6.6.4]{HP}).     

\begin{lemma}\label{L-3.4}
For each $t\geq1$ and $\varepsilon>0$ there exist $N_0, m_0 \in \mathbb{N}$,
$\eps_0 >0$ and a real polynomial $G$ such that
$\Vert v_{p,q} - (1-p)qp\cdot G(pqp) \Vert_t < \varepsilon$ and such that, for each
$N \geq N_0$, if $(P,Q) \in (G(N,k(N))\times G(N,l(N)))_0$ and if 
\begin{equation}\label{F-3.3}
\left|\frac{1}{N}\mathrm{Tr}_N((PQP)^m) - \tau((pqp)^m)\right| < \eps_0 
\quad \mbox{for $1\leq m \leq m_0$}, 
\end{equation} 
then $\Vert V_{P,Q} - (1-P)QP\cdot G(PQP)\Vert_t < \varepsilon$. Here,
$\|\cdot\|_t$ denotes the Schatten $t$-norm with respect to $\tau$ as well as
$N^{-1}\Tr_N$.
\end{lemma}

\begin{proof}
We only sketch the proof since it is essentially similar to that of
\cite[Lemma 2.6]{HP1}. For each small $\alpha,\beta>0$ we estimate 
\begin{align}\label{F-3.4}
&\Vert v_{p,q} - ((1-p)qp)(\sqrt{pqp(p-pqp)}+\alpha1)^{-1}\Vert_t^t \nonumber\\
&\quad\leq 
\frac{1}{2}\left\{\nu((0,\beta)) + \nu((1-\beta,1)) + 
\nu([\beta,1-\beta])\left(\frac{\alpha}{\sqrt{\beta(1-\beta)}+\alpha}\right)^t
\right\}
\end{align}
and
\begin{align}\label{F-3.5}
&\Vert V_{P,Q} - (I-P)QP(\sqrt{PQP(P-PQP)}+\alpha I)^{-1}\Vert_t^t \nonumber\\
&\quad\leq 
\frac{1}{N}\#\left\{i : \lambda_i(PQP) < \beta\right\}
+ \frac{1}{N}\#\left\{i : \lambda_i(PQP) > 1-\beta\right\} \nonumber\\
&\quad\qquad + \frac{k(N)}{N}\left(\frac{\alpha}{\sqrt{\beta(1-\beta)}+\alpha}
\right)^t, 
\end{align} 
where $0 < \lambda_1(PQP) < \cdots < \lambda_{k(N)}(PQP) < 1$ are the eigenvalues 
of $PQP|_{P\bC^N}$ for $(P,Q) \in (G(N,k(N))\times G(N,l(N)))_0$. For any $\eta>0$ let
us choose a $\beta>0$ so that $\nu((0,2\beta)) + \nu((1-2\beta,1)) < \eta^t$. By 
\eqref{F-3.4} we get 
\begin{equation}\label{F-3.6}
\Vert v_{p,q} - ((1-p)qp)(\sqrt{pqp(p-pqp)}+\alpha 1)^{-1}\Vert_t^t 
\leq \frac{\eta^t}{2} + \frac{\tau(E)}{2}\left(
\frac{\alpha}{\sqrt{\beta(1-\beta)}}\right)^t.
\end{equation}
Note that $\nu$ is non-atomic on $(0,1)$ due to the assumption
$\chi_{\mathrm{proj}}(p,q) > -\infty$. Set $\xi_{N,i}:=\min\{x\in[0,1]:\nu((0,x))
=i\tau(E)/k(N)\}$ for $1\le i\le k(N)$; then we get
$$
\tau((pqp)^m)= \lim_{N\rightarrow\infty}\frac{1}{N}
\sum_{i=1}^{k(N)} \big(\xi_{N,i}\big)^m
\quad\mbox{for all $m \in \mathbb{N}$}.
$$
Also choose a constant $C>\sup_{N\ge2}N/k(N)$. By \cite[4.3.4]{HP} there
are $m_0 \in \mathbb{N}$ and $\eps_0>0$ such that, for every $N\in\mathbb{N}$ and for
every $(\lambda_1,\dots,\lambda_{k(N)}) \in (0,1)_<^{k(N)}$, 
\begin{equation*} 
\left|\frac{1}{k(N)}\sum_{i=1}^{k(N)} \lambda_i^m
- \frac{1}{k(N)}\sum_{i=1}^{k(N)}\big(\xi_{N,i}\big)^m\right| < 2C\eps_0
\quad\mbox{for $1\le m\le m_0$}
\end{equation*} 
implies 
\begin{equation}\label{F-3.7}
\frac{1}{k(N)} \sum_{i=1}^{k(N)} \big|\lambda_i - \xi_{N,i}\big|^m
< \beta \eta^t. 
\end{equation} 
Assume \eqref{F-3.7}. Set $i_0 := \#\{ i : \lambda_i < \beta\}$ and
$i_1 := \#\{i : \xi_{N,i}< 2\beta\}$. If $i_1 < i \leq i_0$, then
$\big|\lambda_i - \xi_{N,i}\big| = \xi_{N,i} - \lambda_i \geq \beta$ so that
we get $i_0 < i_1 + k(N)\eta^t$ by \eqref{F-3.7}. Since
$i_1\tau(E)/k(N) \leq \nu((0,2\beta)) < \eta^t$, we get
$i_0 < \tau(E)^{-1}(1+\tau(E))k(N)\eta^t$. If there is no $i_1 < i \leq i_0$, then
$i_0 \leq i_1 <\tau(E)^{-1}k(N)\eta^t$. Therefore,
$\#\{i:\lambda_i<\beta\}<\tau(E)^{-1}(1+\tau(E))k(N)\eta^t$. Similarly, we have
$\#\{i:\lambda_i>1-\beta\}<\tau(E)^{-1}(1+\tau(E))k(N)\eta^t$.
Now, choose an $N_0 \in \mathbb{N}$ so that  
\begin{equation*} 
\left|\frac{1}{N}\sum_{i=1}^{k(N)} \big(\xi_{N,i}\big)^m
- \tau((pqp)^m)\right| < \eps_0
\end{equation*} 
for all $1 \leq m \leq m_0$ and $N \geq N_0$. We then conclude that, for every
$N\ge N_0$, if $(P,Q) \in (G(N,k(N))\times G(N,l(N)))_0$ satisfies \eqref{F-3.3}, then  
\begin{align} 
\#\{ i : \lambda_i(PQP) < \beta \} &< \frac{1+\tau(E)}{\tau(E)}k(N)\eta^t, 
\label{F-3.8}\\ 
\#\{ i : \lambda_i(PQP) > 1-\beta \} &< \frac{1+\tau(E)}{\tau(E)}k(N)\eta^t. 
\label{F-3.9}
\end{align}
Inserting \eqref{F-3.8} and \eqref{F-3.9} in \eqref{F-3.5} we get
\begin{align}
&\Vert V_{P,Q} - (I-P)QP(\sqrt{PQP(P-PQP)}+\alpha I)^{-1}\Vert_t^t \nonumber\\
&\qquad\leq \frac{1+\tau(E)}{\tau(E)}\eta^t + \frac{1}{2}\left(
\frac{\alpha}{\sqrt{\beta(1-\beta)}}\right)^t. \label{F-3.10}
\end{align}
Finally, let $\alpha>0$
be so small as $\alpha/\sqrt{\beta(1-\beta)} < \eta$, and choose a real polynomial
$G(x)$ such that $|G(x) - (\sqrt{x(1-x)}+\alpha)^{-1}| < \eta$ for all $x \in [0,1]$.
Then by \eqref{F-3.6} and \eqref{F-3.10} we obtain
$$ 
\Vert v_{p,q} - (1-p)qp\cdot G(pqp)\Vert_t < 2\eta 
$$
and
$$
\Vert V_{P,Q} - (I-P)QP\cdot G(PQP)\Vert_t 
< \left(\left(\frac{1}{\tau(E)}+{3\over2}\right)^{1/t}+1\right) \eta.
$$
The proof is completed if $\eta>0$ was chosen
so small as $\left((1/\tau(E)+3/2)^{1/t}+1\right)\eta < \varepsilon$. 
\end{proof}  

\medskip\noindent
{\it Proof of Lemma \ref{L-3.3}.}\enspace
Choose $k_1(N),\dots,k_{n'}(N)$ so that $k_i(N)/N \rightarrow \tau(r_i)$ as
$N\rightarrow\infty$, and set
$$
\Phi_N := \Phi_{N,\psi}\times\prod_{i=1}^{n'}\mathrm{id}_{G(N,k_i(N))}
\quad\mbox{on}\ (G(N,k(N))\times G(N,l(N)))_0\times\prod_{i=1}^{n'}G(N,k_i(N))
$$
and $\gamma_N := \gamma_{G(N,k(N))}\otimes\gamma_{G(N,l(N))}\otimes
\bigotimes_{i=1}^{n'}\gamma_{G(N,k_i(N))}$. Let $m \in \mathbb{N}$ and $\varepsilon>0$
be arbitrary. In the following, for brevity we write
$\Gamma_\proj(p,q,r_1,\dots,r_{n'};N,m_0,\eps_0)$ etc.\ without
$k(N),l(N),k_1(N),\dots,k_n(N)$. Thanks to Lemma \ref{L-3.4} together with the expressions of
$q(\psi;p)$ and $Q(\psi;P)$ above, we can choose $N_0, m_0 \in \mathbb{N}$ and
$\varepsilon_0>0$ with $m_0 \geq m$ and $\varepsilon_0 \leq \varepsilon$ such that, for every
$N \geq N_0$, if $(P,Q,R_1,\dots,R_{n'}) \in \Gamma_{\mathrm{proj}}(p,q,r_1,\dots,r_{n'};
N,m_0,\varepsilon_0)$ and $(P,Q)\in(G(N,k(N))\times G(N,l(N)))_0$, then
$\Phi_N(P,Q,R_1,\allowbreak\dots,R_{n'})$ falls into
$\Gamma_{\mathrm{proj}}(p,q(\psi;p),r_1,\dots,r_{n'};\allowbreak N,m,\varepsilon)$.
Via $\Xi_{N,k(N),l(N)}$ in the first two coordinates, Lemma \ref{L-3.1} enables us to
estimate the Radon-Nikodym derivative $d\gamma_N\circ\Phi_N/d\gamma_N$ on a
co-negligible subset of $\Gamma_{\mathrm{proj}}(p,q,r_1,\dots,r_{n'};N,m,\varepsilon)$
from below by the infimum value of 
\begin{align}\label{F-3.11}
&\prod_{1\leq i < j \leq k(N)}
\left(\frac{\psi(\lambda_i(PQP))-\psi(\lambda_j(PQP))}
{\lambda_i(PQP)-\lambda_j(PQP)}\right)^2
\,\prod_{i=1}^{k(N)}\psi'(\lambda_i(PQP)) \notag\\
&\qquad\times \prod_{i=1}^{k(N)} \left(\frac{\psi(\lambda_i(PQP))}{\lambda_i(PQP)}
\right)^{l(N)-k(N)}
\,\prod_{i=1}^{k(N)}\left(\frac{1-\psi(\lambda_i(PQP))}{1-\lambda_i(PQP)}
\right)^{N-k(N)-l(N)} 
\end{align} 
for all $(P,Q) \in (G(N,k(N))\times G(N,l(N)))_0\cap\Gamma_{\mathrm{proj}}
(p,q;N,m_0,\varepsilon_0)$ with the eigenvalue list
$\lambda_1(PQP),\dots,\lambda_{k(N)}(PQP)$ in increasing order. 

Let $\psi^{[1]}(x,y)$ be the so-called divided quotient of $\psi$, i.e., 
\begin{equation*} 
\psi^{[1]}(x,y) := 
\begin{cases} 
\frac{\psi(x)-\psi(y)}{x-y} & (x\neq y), \\
\psi'(x) & (x=y). 
\end{cases}
\end{equation*}
Then, quantity \eqref{F-3.11} is rewritten in the coordinate $(P,Q)$ as 
\begin{align*} 
&\mathrm{det}_{k(N)^2\times k(N)^2}\left[ P\otimes P\cdot
\psi^{[1]}(PQP\otimes P,P\otimes PQP)\cdot P\otimes P \right] \\
&\quad\times \left(\mathrm{det}_{k(N)\times k(N)}[P(PQP)^{-1}\psi(PQP)P]
\right)^{l(N)-k(N)} \\
&\quad\times \left(\mathrm{det}_{k(N)\times k(N)}[P(P-PQP)^{-1}(P-\psi(PQP)P]
\right)^{N-k(N)-l(N)} \\ 
&= \exp\left(\mathrm{Tr}_{k(N)}^{\otimes2}\left(P\otimes P\cdot
\log(\psi^{[1]}(PQP\otimes P, P\otimes PQP))\cdot P\otimes P\right) \right) \\
&\quad\times
\left(\exp\left(\mathrm{Tr}_{k(N)}\left(P\cdot\log\left((PQP)^{-1}\psi(PQP)
\right)\cdot P\right)\right)\right)^{l(N)-k(N)} \\
&\quad\times
\left(\exp\left(\mathrm{Tr}_{k(N)}\left(P\cdot\log\left((P-PQP)^{-1}(P-\psi(PQP))
\right)\cdot P\right)\right)\right)^{N-k(N)-l(N)},
\end{align*}
where $\psi^{[1]}(PQP\otimes P,P\otimes PQP)$ is defined on $P\bC^N\otimes P\bC^N$
while $(PQP)^{-1}\psi(PQP)$ and $(P-PQP)^{-1}(P-\psi(PQP))$ are on $P\bC^N$. Let
$\delta>0$ be arbitrary. Since $\psi$ is $C^1$, $\log\psi^{[1]}(x,y)$ is continuous
on $[0,1]^2$ so that there is a real polynomial $L(x,y)$ on $[0,1]^2$ such that
$\Vert \log\psi^{[1]} - L\Vert_{\infty} < \delta$. If $m'\in\mathbb{N}$ is larger
than the degree of $L$, then we have, for each
$(P,Q) \in \Gamma_{\mathrm{proj}}(p,q;N,m',\varepsilon')$ with an arbitrary
$\varepsilon'>0$,   
\begin{align*} 
&\Big|\frac{1}{N^2}\mathrm{Tr}_{N}^{\otimes2}\left(P\otimes P\cdot
\log\psi^{[1]}(PQP\otimes P, P\otimes PQP)\cdot P\otimes P\right) \\
&\phantom{aaaaaaaaaaaaaaa}-
\tau^{\otimes2}(p\otimes p\cdot \log\psi^{[1]}(pqp\otimes p, p\otimes pqp)
\cdot p\otimes p)\Big| \\
&\quad\leq 2\delta + 
\Big|\frac{1}{N^2}\mathrm{Tr}_{N}^{\otimes2}(P\otimes P\cdot
L(PQP\otimes P, P\otimes PQP)\cdot P\otimes P) \\
&\phantom{aaaaaaaaaaaaaaaaaaa}- 
\tau^{\otimes2}(p\otimes p\cdot L(pqp\otimes p, p\otimes pqp)\cdot
p\otimes p)\Big| \\
&\quad\leq 
2\delta+C\varepsilon     
\end{align*}
with $C>0$ depending only on $L$ (hence on $\delta$). Therefore, for  each $\eta>0$
there are $m_1 \in \mathbb{N}$ and $\varepsilon_1>0$ such that
\begin{align} 
\exp&\left(\mathrm{Tr}_{k(N)}^{\otimes2}\left(P\otimes P\cdot
\log\left(\psi^{[1]}(PQP\otimes P, P\otimes PQP)\right)\cdot
P\otimes P\right)\right) \notag\\
&\geq\exp\left(N^2\left\{\tau^{\otimes2}(p\otimes p\cdot
\log\psi^{[1]}(pqp\otimes p, p\otimes pqp)\cdot
p\otimes p)-\eta\right\}\right) \label{F-3.12}
\end{align}   
for all $(P,Q) \in \Gamma_{\mathrm{proj}}(p,q;N,m',\varepsilon')$
as long as $m' \geq m_1$ and $0<\varepsilon'\leq\varepsilon_1$. Since $x^{-1}\psi(x)$
and $(1-x)^{-1}(1-\psi(x))$ are both bounded below above $0$ on $[0,1]$ due to the
assumption on $\psi$, the same argument works for the other two terms
\begin{gather*} 
\exp\left(\mathrm{Tr}_{k(N)}\left(P\cdot\log\left((PQP)^{-1}\psi(PQP)
\right)\cdot P\right)\right), \\
\exp\left(\mathrm{Tr}_{k(N)}\left(P\cdot\log\left((P-PQP)^{-1}(P-\psi(PQP))
\right)\cdot P\right)\right).
\end{gather*}  
Therefore, for each $\eta>0$ there are $m_2 \in \mathbb{N}$ and $\varepsilon_2>0$
such that 
\begin{align} 
\exp&\left(\mathrm{Tr}_{k(N)}\left(P\cdot\log\left((PQP)^{-1}\psi(PQP)
\right)\cdot P\right)\right) \notag\\
&\geq
\exp\left(N\left\{\tau(p\cdot\log((pqp)^{-1}\psi(pqp))\cdot p)-\eta\right\}\right),
\label{F-3.13}\\
\exp&\left(\mathrm{Tr}_{k(N)}\left(P\cdot\log\left((P-PQP)^{-1}(P-\psi(PQP))
\right)\cdot P\right)\right) \notag\\
&\geq\exp\left(N\left\{\tau(p\cdot\log((p-pqp)^{-1}(p-\psi(pqp)))\cdot p)
-\eta\right\}\right) \label{F-3.14}
\end{align}
for all $(P,Q) \in \Gamma_{\mathrm{proj}}(p,q;N,m',\varepsilon')$ as long as
$m'\geq m_2$ and $0<\varepsilon'\leq\varepsilon_2$. Hence, for every $N \geq N_0$,
$m'\geq \max\{m_0,m_1,m_2\}$ and
$0<\varepsilon'<\min\{\varepsilon_0,\varepsilon_1,\varepsilon_2\}$, we have
\allowdisplaybreaks{
\begin{align*} 
&\frac{1}{N^2}\log\gamma_N\bigl(\Gamma_{\mathrm{proj}}
(p,q(\psi;p),r_1,\dots,r_{n'};N,m,\varepsilon)\bigr) \\
&\qquad\geq 
\frac{1}{N^2}\log\gamma_N\bigl(\Phi_N\bigl(\Gamma_{\mathrm{proj}}
(p,q,r_1,\dots,r_{n'};N,m',\varepsilon')\bigr)\bigr) \\
&\qquad\geq 
\frac{1}{N^2}\log\gamma_N\bigl(\Gamma_{\mathrm{proj}}
(p,q,r_1,\dots,r_{n'};N,m',\varepsilon')\bigr) \\
&\qquad\quad+\tau^{\otimes2}(p\otimes p\cdot\log\psi^{[1]}(pqp\otimes p,p\otimes pqp)
\cdot p\otimes p) \\
&\qquad\quad+ \left(\frac{l(N)}{N}-\frac{k(N)}{N}\right)
\tau(p\cdot\log((pqp)^{-1}\psi(pqp))\cdot p) \\
&\qquad\quad+ \left(1-\frac{l(N)}{N}-\frac{k(N)}{N}\right)
\tau(p\cdot\log((p-pqp)^{-1}(p-\psi(pqp)))\cdot p)-3\eta \\
&\qquad= 
\frac{1}{N^2}\log\gamma_N\bigl(\Gamma_{\mathrm{proj}}
(p,q,r_1,\dots,r_{n'};N,m',\varepsilon')\bigr) \\
&\qquad\quad+ \frac{1}{4}\iint_{(0,1)^2} \log\left|\frac{\psi(x)-\psi(y)}{x-y}\right|
\,d\nu(x)\,d\nu(y) \\
&\qquad\quad+ \frac{1}{2}\left(\frac{l(N)}{N}-\frac{k(N)}{N}\right)
\int_{(0,1)} \log\frac{\psi(x)}{x}\,d\nu(x) \\
&\qquad\quad+ \frac{1}{2}\left(1-\frac{l(N)}{N}-\frac{k(N)}{N}\right)
\int_{(0,1)} \log\frac{1-\psi(x)}{1-x}\,d\nu(x) - 3\eta.  
\end{align*} 
}Take the $\limsup$ as $N\rightarrow\infty$ and the limit as $m\rightarrow\infty$,
$\varepsilon\searrow0$ in the above inequality. Since $\eta>0$ is arbitrary, we get
\begin{align*} 
&\chi_{\mathrm{proj}}(p,q(\psi;p),r_1,\dots,r_{n'}) \\
&\quad\geq \chi_{\mathrm{proj}}(p,q,r_1,\dots,r_{n'})
+ \frac{1}{4}\iint_{(0,1)^2} \log\left|\frac{\psi(x)-\psi(y)}{x-y}\right|
\,d\nu(x)d\nu(y) \\
&\qquad+ {\tau(q)-\tau(p)\over2}
\int_{(0,1)} \log\frac{\psi(x)}{x}\,d\nu(x)
+ {1-\tau(q)-\tau(p)\over2}
\int_{(0,1)} \log\frac{1-\psi(x)}{1-x}\,d\nu(x) \\
&\quad =\chi_\proj(p,q,r_1,\dots,r_{n'})
+\chi_\proj(p,q(\psi;p))-\chi_\proj(p,q) 
\end{align*}
thanks to Proposition \ref{P-2.1}. The reverse inequality can be shown as well if we 
replace the inequalities \eqref{F-3.12}--\eqref{F-3.14} by their reversed versions.
\qed

\bigskip
For the second step we present two more technical lemmas. The proof of the next lemma
should be compared with that of \cite[Lemma 4.1]{V4}.

\begin{lemma}\label{L-3.5}
Let $\mu$ be a measure on $[0,1]$ with no atom at $0,1$, and assume the conditions
\begin{gather} 
\iint_{(0,1)^2} \log|x-y|\,d\mu(x)\,d\mu(y) > -\infty, \label{F-3.15}\\
\int_{(0,1)} \log x\, d\mu(x) > -\infty, \label{F-3.16}\\ 
\int_{(0,1)} \log(1-x)\,d\mu(x) > -\infty. \label{F-3.17}
\end{gather}
If $\psi$ is a continuous increasing function from $[0,1]$ onto itself with
$\psi(0) = 0$, $\psi(1) = 1$, then there exists a sequence of
$C^{\infty}$-diffeomorphisms $\psi_j$ from $[0,1]$ onto itself with $\psi_j(0) = 0$,
$\psi_j(1) = 1$ such that 
\begin{itemize} 
\item[(i)] $\psi'_j(x)\geq1/j$ for all $j\in\bN$ and $x \in [0,1]$,
\item[(ii)] $\psi_j \longrightarrow \psi$ uniformly on $[0,1]$,
\item[(iii)] $\displaystyle{\lim_{j\rightarrow\infty}\iint_{(0,1)^2} \log|x-y|
\,d(\psi_j{}_*\mu)(x)\,d(\psi_j{}_*\mu)(y) = \iint_{(0,1)^2} \log|x-y|
\,d(\psi_*\mu)(x)\,d(\psi_*\mu)(y)}$,
\item[(iv)] $\displaystyle{\lim_{j\rightarrow\infty}\int_{(0,1)} \log x
\,d(\psi_j{}_*\mu)(x) = \int_{(0,1)} \log x\,d(\psi_*\mu)(x)}$,
\item[(v)] $\displaystyle{\lim_{j\rightarrow\infty}\int_{(0,1)} \log (1-x)
\,d(\psi_j{}_*\mu)(x) = \int_{(0,1)} \log (1-x)\,d(\psi_*\mu)(x)}$,
\end{itemize}
where $\psi_*\mu$ is the image measure of $\mu$ by $\psi$.
Furthermore, when conditions \eqref{F-3.16} and/or \eqref{F-3.17} for $\mu$ are
dropped, the conclusion holds without {\rm(iv)} and/or {\rm(v)} correspondingly.
\end{lemma}

\begin{proof}
Extend $\psi$ to a continuous increasing function on the whole $\bR$ periodically,
namely, $\psi(x+m) = \psi(x)+m$ for $x \in [0,1]$ and $m\in\mathbb{Z}$.
For each $j \in \mathbb{N}$, by \eqref{F-3.15}--\eqref{F-3.17} one can choose a
$\delta_j \in(0,1/j]$ such that 
\begin{align} 
\iint_{\{ (x,y) \in (0,1)^2 : |x-y|<\delta_j\}} \log|x-y|\,d\mu(x)\,d\mu(y)
&\geq -1/j, \label{F-3.18}\\
\iint_{\{ (x,y) \in (0,1)^2 : |x-y|<\delta_j\}} d\mu(x)\,d\mu(y)
&\leq 1/(j\log j), \label{F-3.19}\\
\int_{(0,\delta_j)} \log x\,d\mu(x) &\geq -1/j, \label{F-3.20}\\
\mu((0,\delta_j)) &\leq 1/(j\log j), \label{F-3.21}\\
\int_{(1-\delta_j,1)} \log(1-x)\,d\mu(x) &\geq -1/j, \label{F-3.22}\\
\mu((1-\delta_j,1)) &\leq 1/(j\log j). \label{F-3.23} 
\end{align}
For each $j$ we choose a $C^{\infty}$-function $\phi_j \geq 0$ supported in
$[-1/j,1/j]$ with $\int \phi_j(x)\,dx = 1$ such that 
$|(\psi*\phi_j)(x) - \psi(x)| \leq \delta_j/2j$ for all $x \in [0,1]$, and define
$$
\psi_j(x) := {x\over j} + \biggl(1-{1\over j}\biggr)
((\psi*\phi_j)(x) - (\psi*\phi_j)(0))
\quad\mbox{for $x\in[0,1]$}.
$$
Then one can immediately see that $\psi_j$ is $C^\infty$, $\psi_j(0)=0$, $\psi_j(1)=1$
and (i), (ii) are satisfied.

For $x,y \in [0,1]$ with $|x-y|\geq\delta_j$ notice that
\begin{align*} 
&|\psi_j(x) - \psi_j(y)| \\
&\quad={|x-y|\over j}+\biggl(1-\frac{1}{j}\biggr)
|(\psi*\phi_j)(x)-(\psi*\phi_j)(y)| \\
&\quad\geq {|x-y|\over j}+ \biggl(1-\frac{1}{j}\biggr)
\big\{|\psi(x)-\psi(y)| - |(\psi*\phi_j)(x)-\psi(x)|-|\psi(y)-(\psi*\phi_j)(y)|
\big\} \\
&\quad\ge \biggl(1-\frac{1}{j}\biggr)|\psi(x)-\psi(y)|,
\end{align*}  
and in particular
\begin{align*}
\psi_j(x)\ge&\biggl(1-\frac{1}{j}\biggr)\psi(x)
\quad\qquad\ \,\mbox{for $x \in [\delta_j,1)$}, \\
1-\psi_j(x)\ge&\biggl(1-\frac{1}{j}\biggr)(1-\psi(x))
\quad\mbox{for $x \in (0,1-\delta_j]$}.
\end{align*}
Hence we have by \eqref{F-3.18} and \eqref{F-3.19}
\begin{align*} 
&\iint_{(0,1)^2} \log|x-y|\,d(\psi_j{}_*\mu)(x)d(\psi_j{}_*\mu)(y) \\
&\qquad\geq 
\iint_{\{(x,y) \in (0,1)^2 : |x-y|<\delta_j\}} \log{|x-y|\over j}\,d\mu(x)\,d\mu(y) \\
&\qquad\quad+ \iint_{\{(x,y) \in (0,1)^2 : |x-y|\geq\delta_j\}}\log
\biggl(\biggl(1-\frac{1}{j}\biggr)|\psi(x)-\psi(y)|\biggr)\,d\mu(x)\,d\mu(y) \\
&\qquad\ge -\frac{2}{j} +
\log\biggl(1-\frac{1}{j}\biggr)+\iint_{(0,1)^2}\log|x-y|
\,d(\psi_*\mu)(x)\,d(\psi_*\mu)(y), \\
\end{align*}
and also we have by \eqref{F-3.20}--\eqref{F-3.23}
\begin{align*} 
&\int_{(0,1)} \log x\,d(\psi_j{}_*\mu)(x) \\
&\qquad\geq \int_{(0,\delta_j)} \log{x\over j}\,d\mu(x) 
+ \int_{[\delta_j,1)} \log\biggl(\biggl(1-\frac{1}{j}\biggr)\psi(x)\biggr)\,d\mu(x) \\
&\qquad\ge -\frac{2}{j} + \log\biggl(1-\frac{1}{j}\biggr)
+ \int_{(0,1)} \log \psi(x)\,d\mu(x),
\end{align*}
\begin{align*}
&\int_{(0,1)} \log(1-x)\,d(\psi_j{}_*\mu)(x) \\
&\qquad\geq
\int_{(1-\delta_j,1)} \log{1-x\over j}\,d\mu(x) +
\int_{(0,1-\delta_j]} \log\biggl(\biggl(1-\frac{1}{j}\biggr)(1-\psi(x))\biggr)
\,d\mu(x) \\
&\qquad\ge -\frac{2}{j} + \log\biggl(1-\frac{1}{j}\biggr)
+ \int_{(0,1)} \log(1-\psi(x))\,d\mu(x).
\end{align*}
Therefore,
\begin{align*}
&\liminf_{j\rightarrow\infty}\iint_{(0,1)^2} \log|x-y|
\,d(\psi_j{}_*\mu)(x)\,d(\psi_j{}_*\mu)(y) \\
&\qquad\geq \iint_{(0,1)^2} \log|x-y|\,d(\psi_*\mu)(x)\,d(\psi_*\mu)(y),
\end{align*}
\begin{align*}
\liminf_{j\rightarrow\infty} \int_{(0,1)} \log x\,d(\psi_j{}_*\mu)(x)
&\geq \int_{(0,1)} \log x\, d(\psi_*\mu)(x), \\
\liminf_{j\rightarrow\infty} \int_{(0,1)} \log(1-x)\,d(\psi_j{}_*\mu)(x)
&\geq \int_{(0,1)} \log(1-x)\, d(\psi_*\mu)(x).
\end{align*}
On the other hand, Fatou's lemma says that the reverse inequalities of these three
with $\limsup$ in place of $\liminf$ actually hold true. Hence we have (iii)--(v).
Finally, the above proof shows the last statement as well.  
\end{proof}

\begin{lemma}\label{L-3.6}
Let $\mu$ be a measure on $[0,1]$ with no atom at $0,1$, and $\psi$ be a continuous
increasing function from $[0,1]$ into itself. Assume that $\mu$ satisfies conditions
\eqref{F-3.15}--\eqref{F-3.17} in Lemma \ref{L-3.5} and also $\psi_*\mu$ does
\eqref{F-3.16} and \eqref{F-3.17}. Then, there exists a sequence of
continuous increasing functions $\psi_m$ from $[0,1]$ onto itself with
$\psi_m(0) = 0$, $\psi_m(1) = 1$ such that
\begin{itemize}
\item[(i)] $\int_{(0,1)}|\psi_m(x)-\psi(x)|^2\,d\mu(x) \longrightarrow 0$,
\item[(ii)] $\displaystyle{\lim_{m\rightarrow\infty}\iint_{(0,1)^2} \log|x-y|
\,d(\psi_m{}_*\mu)(x)\,d(\psi_m{}_*\mu)(y) = \iint_{(0,1)^2} \log|x-y|
\,d(\psi_*\mu)(x)\,d(\psi_*\mu)(y)}$,
\item[(iii)] $\displaystyle{\lim_{m\rightarrow\infty}\int_{(0,1)} \log x
\,d(\psi_m{}_*\mu)(x) = \int_{(0,1)} \log x\,d(\psi_*\mu)(x)}$,
\item[(iv)] $\displaystyle{\lim_{m\rightarrow\infty}\int_{(0,1)} \log (1-x)
\,d(\psi_m{}_*\mu)(x) = \int_{(0,1)} \log (1-x)\,d(\psi_*\mu)(x)}$.
\end{itemize}
Furthermore, when conditions \eqref{F-3.16} and/or \eqref{F-3.17} for $\mu$ and
$\psi{}_*\mu$ are dropped, the conclusion holds without {\rm(iii)} and/or {\rm(iv)}
correspondingly.
\end{lemma}

\begin{proof}
We assume that both $\psi(0)>0$ and $\psi(1)<1$; the other cases can be handled
easier. Condition \eqref{F-3.15} implies
\begin{align}
(-\log m)\nu((0,1/m))^2&\ge\iint_{(0,1/m)^2}\log|x-y|\,d\nu(x)\,d\nu(y)
\longrightarrow0, \label{F-3.24}\\
(-\log m)\nu((1-1/m,1))^2&\ge\iint_{(1-1/m,1)^2}\log|x-y|\,d\nu(x)\,d\nu(y)
\longrightarrow0 \label{F-3.25}
\end{align}
as $m\to\infty$. On the other hand, \eqref{F-3.16} and \eqref{F-3.17} imply
\begin{align}
(-\log m)\nu((0,1/m)) &\geq \int_{(0,1/m)} \log x\,d\mu(x) \longrightarrow 0,
\label{F-3.26}\\
(-\log m)\nu((1-1/m,1)) &\geq \int_{(1-1/m,1)}\log(1-x)\,d\mu(x)\longrightarrow 0,
\label{F-3.27}
\end{align}
respectively. For each $m\ge2$ define a function $\psi_m$ on $[0,1]$ by 
\begin{equation*}
\psi_m(x) := 
\begin{cases} 
mx\psi(x) & (0\leq x < 1/m), \\
\psi(x) & (1/m \leq x \leq 1-1/m), \\
1-m(1-x)(1-\psi(x)) & (1-1/m < x \leq 1), 
\end{cases} 
\end{equation*}
which is clearly continuous and increasing with $\psi_m(0)=0$, $\psi_m(1)=1$. Then
(i) immediately follows. It is easy to check the following: 
$$
|\psi_m(x) - \psi_m(y)| \geq \begin{cases}
 m \psi(0) |x-y| & \text{for $x,y\in(0,1/m)$}, \\
m (1-\psi(1))|x-y| & \text{for $x,y\in(1-1/m,1)$}, \\
|\psi(x)-\psi(y)| & \text{for other $x,y\in(0,1)$}.
\end{cases}
$$
Hence we have
\begin{align*} 
&\iint_{(0,1)^2} \log|x-y|\,d(\psi_m{}_*\mu)(x)\,d(\psi_m{}_*\mu)(y) \\
&\quad\geq  
\iint_{(0,1/m)^2} \log (m\psi(0)|x-y|)\,d\mu(x)\,d\mu(y) \\
&\qquad+ \iint_{(1-1/m,1)^2} \log(m(1-\psi(1))|x-y|)\,d\mu(x)\,d\mu(y) \\
&\qquad+ \iint_{(0,1)^2} \log|\psi(x)-\psi(y)|\,d\mu(x)\,d\mu(y) \\
&\quad=(\log m + \log \psi(0))\mu((0,1/m))^2
+ \iint_{(0,1/m)^2} \log|x-y|\,d\mu(x)\,d\mu(y) \\
&\qquad+(\log m + \log(1-\psi(1)))\mu((1-1/m,1))^2
+ \iint_{(1-1/m,1)^2} \log|x-y|\,d\mu(x)\,d\mu(y) \\
&\qquad+ \iint_{(0,1)^2} \log|\psi(x)-\psi(y)|\,d\mu(x)\,d\mu(y) \\
&\quad\longrightarrow \iint_{(0,1)^2} \log|\psi(x)-\psi(y)|\,d\mu(x)\,d\mu(y)
\end{align*}
as $m\to\infty$ by \eqref{F-3.24}, \eqref{F-3.25} and \eqref{F-3.15}. Therefore,
\begin{align*}
&\liminf_{m\rightarrow\infty}\iint_{(0,1)^2}\log|x-y|
\,d(\psi_m{}_*\mu)(x)\,d(\psi_m{}_*\mu)(y) \\
&\qquad\geq \iint_{(0,1)^2}\log|x-y|\,d(\psi_m{}_*\mu)(x)\,d(\psi_m{}_*\mu)(y). 
\end{align*}
This together with Fatou's lemma implies (ii). On the other hand, by \eqref{F-3.16}
for $\mu$ and $\psi{}_*\mu$ we have
\begin{align*} 
&\int_{(0,1/m)} \log(mx\psi(x))\,d\mu(x) \\
&\qquad=(\log m) \nu((0,1/m)) + \int_{(0,1/m)} \log x\,d\mu(x) 
+ \int_{(0,1/m)} \log \psi(x)\,d\mu(x)\longrightarrow 0
\end{align*}
thanks to \eqref{F-3.26}.
Furthermore,
\begin{align*}
0&\ge\int_{(1-1/m,1)} \log(1-m(1-x)(1-\psi(x)))\,d\mu(x) \\
&\geq \log \psi(1-1/m) \cdot \mu((1-1/m,1)) \longrightarrow 0.
\end{align*}
These imply (iii). Similarly, (iv) follows from \eqref{F-3.17} for $\mu$ and
$\psi{}_*\mu$ thanks to \eqref{F-3.27}.
\end{proof}
 
We are now in the final position to prove Theorem \ref{T-3.2} in full generality.

\medskip\noindent
{\it Proof of Theorem \ref{T-3.2}.}\enspace
As mentioned before we may assume $n=1$, and write $p=p_1$, $q=q_1$ and $\psi=\psi_1$.
We may further assume that $\chi_{\mathrm{proj}}(p,q(\psi;p)) > -\infty$ as well as
$\chi_\proj(p,q)>-\infty$; otherwise, both sides of the inequality are $-\infty$
thanks to Proposition \ref{P-1.2}\,(ii). By Proposition \ref{P-2.1} both $\nu$ and
$\psi{}_*\nu$ satisfy condition \eqref{F-3.15}; moreover they satisfy \eqref{F-3.16}
unless $\tau(p)=\tau(q)$ and also \eqref{F-3.17} unless $\tau(p)=\tau(\1-q)$. In each
case where those equalities of traces occur or not, we choose a sequence $\psi_m$
correspondingly as mentioned in Lemma \ref{L-3.6}. Since
\begin{equation*} 
\Vert p\psi_m(pqp)p-p\psi(pqp)p\Vert_2^2
= \int_{(0,1)} |\psi_m(x) - \psi(x)|^2\,d\nu(x) \longrightarrow 0,
\end{equation*}
we get $p\psi_m(pqp)p\rightarrow p\psi(pqp)p$ strongly so that
$q(\psi_m;p)\rightarrow q(\psi;p)$ strongly as $m\rightarrow\infty$ due to the
definition of $q(\psi;p)$. By Propositions \ref{P-1.2}\,(iii) and \ref{P-2.1} we see
that it suffices to prove the inequality in the case where $\psi(0)=0$ and $\psi(1)=1$.
The same argument using Lemma \ref{L-3.5} in turn enables us to reduce the proof to
Lemma \ref{L-3.3}, and the proof of the inequality is completed.

To prove the equality of the last statement, let $\psi$ be strictly increasing on
$(0,1)$ and define $\tilde \psi$ on $[0,1]$ by
$$
\tilde \psi(x):=\begin{cases}
0 & \text{($0\le x\le \psi(0+)$)}, \\
\psi^{-1}(x) & \text{($\psi(0+)<x<\psi(1-)$)}, \\
1 & \text{($\psi(1-)\le x\le1$)}.
\end{cases}
$$
Furthermore, set $\tilde q:=q(\psi;p)$ and $\tilde\nu:=\psi{}_*\nu$. Then it is clear
that $\tilde\nu$ is the measure corresponding to the pair $(p,\tilde q)$ so that
$\nu=\tilde \psi{}_*\tilde\nu$ and $q=\tilde q(\tilde \psi;p)$. Hence the inequality
established above can be applied to $(p,\tilde q)$ and $\tilde \psi$ too, and we have
the reversed inequality as well.\qed

\section{Additivity and freeness}
\setcounter{equation}{0}

In this section, we prove the next additivity theorem asserting that the pair-block
freeness of projections is characterized by the additivity of their free entropy. For
the projection version of free entropy we have no counterpart of the so-called
infinitesimal change of variable formula in \cite[Proposition 1.3]{V4}, and
hence we need to find another route to prove that the additivity implies the freeness.   

\begin{thm}\label{T-4.1}
Let $p_1,q_1,\dots,p_n,q_n,r_1,\dots,r_{n'}$ be projections in $(\cM,\tau)$.
\begin{itemize}
\item[(1)] If $\{p_1,q_1\}$, $\dots$, $\{p_n,q_n\}$, $\{r_1\}$,
$\dots$, $\{r_{n'}\}$ are free, then
$$
\chi_\proj(p_1,q_1,\dots,p_n,q_n,r_1,\dots,r_{n'}) 
=\chi_\proj(p_1,q_1)+\dots+\chi_\proj(p_n,q_n).
$$
\item[(2)] Conversely, if $\chi_\proj(p_i,q_i)>-\infty$ for $1\le i\le n$ and
equality holds in {\rm(1)}, then $\{p_1,q_1\}$, $\dots$, $\{p_n,q_n\}$, $\{r_1\}$,
$\dots$, $\{r_{n'}\}$ are free.
\item[(3)] In particular, $\chi_\proj(p_1,\dots,p_n)=0$ if and only if $p_1,\dots,p_n$
are free.
\end{itemize}
\end{thm}

\begin{proof}
(1)\enspace It suffices to prove the following two assertions:
\begin{itemize}
\item[(a)] If $\{p,q\}$ and $\{p_1,\dots,p_n\}$ are free, then
$$
\chi_\proj(p,q,p_1,\dots,p_n)=\chi_\proj(p,q)+\chi_\proj(p_1,\dots,p_n).
$$
\item[(b)] If $\{p\}$ and $\{p_1,\dots,p_n\}$ are free, then
$$
\chi_\proj(p,p_1,\dots,p_n)=\chi_\proj(p_1,\dots,p_n).
$$
\end{itemize}
The proofs of these being same, we give only that of (a), which is essentially same as
in \cite{V2,V-IMRN} (see also \cite[pp.\ 269--272]{HP}).

To prove (a), we may assume that $\chi_\proj(p,q)>-\infty$ and
$\chi_\proj(p_1,\dots,p_n)>-\infty$. Choose $k(N),l(N),k_i(N)\in\{0,1,\dots,N\}$ for
$N\in\bN$ and $1\le i\le n$ such that $k(N)/N\to\tau(p)$, $l(N)/N\to\tau(q)$ and
$k_i(N)/N\to\tau(p_i)$ as $N\to\infty$. For each $m\in\bN$ and $\eps>0$ we set
\begin{align*}
\Omega_N(m,\eps)&:=\Gamma_\proj(p,q;k(N),l(N);N,m,\eps) \\
&\qquad\times\Gamma_\proj(p_1,\dots,p_n;k_1(N),\dots,k_n(N);N,m,\eps), \\
\Theta_N(m,\eps)&:=\Gamma_\proj(p,q,p_1,\dots,p_n;k(N),l(N),
k_1(N),\dots,k_n(N);N,m,\eps).
\end{align*}
For given $m\in\bN$ and $\eps>0$ one can show as in \cite[6.4.3]{HP} that
there exists an $\eps_1>0$ such that
$$
\lim_{N\to\infty}{\gamma_N(\Omega_N(m,\eps_1)\cap\Theta_N(m,\eps))
\over\gamma_N(\Omega_N(m,\eps_1))}=1,
$$
where 
$\gamma_N:=\gamma_{G(N,k(N))}\otimes\gamma_{G(N,l(N))}\otimes\gamma_{\vec k(N)}$ and
$\gamma_{\vec k(N)}:=\bigotimes_{i=1}^n\gamma_{G(N,k_i(N))}$. Hence we have
\begin{align*}
&\limsup_{N\to\infty}{1\over N^2}\log\gamma_N(\Theta_N(m,\eps)) \\
&\qquad\ge\limsup_{N\to\infty}{1\over N^2}\log\gamma_N(\Omega_N(m,\eps_1)) \\
&\qquad=\lim_{N\to\infty}{1\over N^2}
\log\bigl(\gamma_{G(N,k(N))}\otimes\gamma_{G(N,l(N))}\bigr)
\bigl(\Gamma_\proj(p,q;k(N),l(N);N,m,\eps_1)\bigr) \\
&\qquad\qquad
+\limsup_{N\to\infty}{1\over N^2}
\log\gamma_{\vec k(N)}
\bigl(\Gamma_\proj(p_1,\dots,p_n;k_1(N),\dots,k_n(N);N,m,\eps_1)\bigr) \\
&\qquad\ge\chi_\proj(p,q)+\chi_\proj(p_1,\dots,p_n).
\end{align*}
The above equality is due to \cite[Proposition 3.3]{HP1}. Therefore,
$$
\chi_\proj(p,q,p_1,\dots,p_n)\ge\chi_\proj(p,q)+\chi_\proj(p_1,\dots,p_n),
$$
and the reverse inequality is Proposition 1.2\,(ii).

\medskip
(3) will be proven in Corollary \ref{C-5.5} of the next section as a
consequence of a transportation cost inequality for projection multi-variables.

\medskip
(2)\enspace We may assume that $p_1,q_1,\dots,p_n,q_n$ are all non-zero. For
$1\le i\le n$ let $\nu_i$ be the measure on $(0,1)$ corresponding to the pair
$(p_i,q_i)$ (see \S2). For each $i$, since $\nu_i$ is non-atomic by the assumption
$\chi_\proj(p_i,q_i)>-\infty$, one can choose a continuous increasing function
$\psi_i$ from $(0,1)$ into itself such that $\psi_i{}_*\nu_i$ is equal to
\eqref{F-2.4} with $\alpha=\tau(p_i)$, $\beta=\tau(q_i)$. Consider $q_i(\psi_i;p_i)$
constructed from $(p_i,q_i)$ and $\psi_i$ (see \S3). Since $\psi_i{}_*\nu_i$
corresponds to the pair $(p_i,q_i(\psi_i;p_i))$, we get
$\chi_\proj(p_i,q_i(\psi_i;p_i))=0$. Therefore, by Theorem \ref{T-3.2} and the
additivity assumption, we have
\begin{align*}
&\chi_\proj(p_1,q_1(\psi_1;p_1),\dots,p_n,q_n(\psi_n;p_n),r_1,\dots,r_{n'}) \\
&\qquad\ge\chi_\proj(p_1,q_1,\dots,p_n,q_n,r_1,\dots,r_{n'})
-\sum_{i=1}^n\chi_\proj(p_i,q_i)=0.
\end{align*}
This implies by (3) that
$p_1,q_1(\psi_1;p_1),\dots,p_n,q_n(\psi_n;p_n),r_1,\dots,r_{n'}$ are free. Since
$\nu_i$ and $\psi_i{}_*\nu_i$ are non-atomic, it is plain to see that
$\{p_i,q_i\}''=\{p_i,q_i(\psi_i;p_i)\}''$ for $1\le i\le n$. Hence the freeness of
$\{p_1,q_1\},\dots,\{p_n,q_n\},\{r_1\},\dots,\{r_{n'}\}$ is obtained.
\end{proof}

\section{Asymptotic freeness and free transportation cost inequality}
\setcounter{equation}{0}

The aim of this section is to prove a transportation inequality for tracial
distributions of projection multi-variables. To do so, we first present an asymptotic
freeness result for random projection matrices generalizing Voiculescu's result in
\cite{V0}.

\subsection{Asymptotic freeness for random projection matrices}
Let $\bigl(\{P(s,N),Q(s,N)\}\bigr)_{s\in S}$ be an independent family of pairs of
$N\times N$ random projection matrices, and let $k(s,N)$, $l(s,N)$, $n_{11}(s,N)$,
$n_{10}(s,N)$, $n_{01}(s,N)$ and $n_{00}(s,N)$ denote the ranks of $P(s,N)$, $Q(s,N)$,
$P(s,N)\wedge Q(s,N)$, $P(s,N)\wedge Q(s,N)^\perp$, $P(s,N)^\perp\wedge Q(s,N)$,
$P(s,N)^\perp\wedge Q(s,N)^\perp$, respectively. For each $s \in S$ we assume the
following: 
\begin{itemize} 
\item[(1)] $k(s,N)$, $l(s,N)$ and $n_{ij}(s,N)$'s are constant almost surely and
$k(s,N)/N$, $l(s,N)/N$ and $n_{ij}(s,N)/N$ converge as $N \rightarrow \infty$. 
\item[(2)] The joint distribution of $(P(s,N),Q(s,N))$ is invariant under unitary
conjugation $(P,Q) \mapsto (UPU^*,UQU^*)$ for $U\in\U(N)$.
\item[(3)] For each $s\in S$ the distribution measure of $P(s,N)Q(s,N)P(s,N)$ with
respect to $N^{-1}\mathrm{Tr}_N$ converges almost surely to a (non-random) measure on
$[0,1]$ as $N\rightarrow\infty$.   
\end{itemize} 
Let $(R(s',N))_{s'\in S'}$ be an independent family of $N\times N$ random projection
matrices, also independent of $\big(\{P(s,N),Q(s,N)\}\big)_{s\in S}$, and assume
that each $R(s',N)$ is distributed under the Haar probability measure on
$G(N,k(s',N))$ with $0\le k(s',N) \leq N$ such that $k(s',N)/N$ converges.
Finally, let $(D(t,N))_{t\in T}$ be a family of $N\times N$ constant matrices such
that $\sup_N \Vert D(t,N)\Vert_{\infty} < +\infty$ for each $t \in T$ and
$(D(t,N),D(t,N)^*)_{t\in T}$ has the limit distribution. In this setup, we have the
following asymptotic freeness result for random projection matrices generalizing
\cite[Theorem 3.11]{V0}.

\begin{thm}\label{T-5.1} With the above notations and assumptions the family 
\begin{equation*} 
\Big(\big(\{P(s,N),Q(s,N)\}\big)_{s\in S},\,\big(R(s',N)\big)_{s' \in S'},
\,\bigl\{D(t,N), D(t,N)^* : t \in T\bigr\}\Big)
\end{equation*} 
is asymptotically free almost surely as $N\rightarrow\infty$.  
\end{thm}

\begin{proof}
Set $n(s,N) := \bigl(N-\sum_{i,j=0}^1 n_{ij}(s,N)\bigr)/2$. By assumption (1),
$n(s,N)$ is constant almost surely and $n(s,N)/N$ converges as $N\to\infty$. As before,
the sine-cosine decomposition of two projections enables us to represent 
\begin{align*} 
P(s,N) &= U(s,N)\left(\begin{bmatrix} I & 0 \\ 0 & 0 \end{bmatrix}
\oplus I \oplus I \oplus 0 \oplus 0\right)U(s,N)^*, \\
Q(s,N) &= U(s,N)\left(\begin{bmatrix} X & \sqrt{X(I-X)} \\ \sqrt{X(I-X)} & I-X
\end{bmatrix}\oplus I\oplus0\oplus I\oplus 0\right)U(s,N)^* 
\end{align*} 
in $\mathbb{C}^N = (\mathbb{C}^{n(s,N)}\otimes\mathbb{C}^2)
\oplus\mathbb{C}^{n_{11}(s,N)}\oplus\mathbb{C}^{n_{10}(s,N)}
\oplus\mathbb{C}^{n_{01}(s,N)}\oplus\mathbb{C}^{n_{00}(s,N)}$, where $U(s,N)$ is a
random unitary matrix and $X=X(s,N)$ is a diagonal matrix whose diagonal entries are
$0\le x_1(s,N)\leq x_2(s,N)\leq\cdots\leq x_{n(s,N)}(s,N)\le1$. Also, we can represent
\begin{equation*} 
R(s',N) = U(s',N)P_{k(s',N)}U(s',N)^*
\end{equation*} 
for each $s'\in S'$, where $U(s',N)$ is a unitary random matrix and $P_{k(s',N)}$ the
diagonal matrix whose first $k(s',N)$ entries are $1$ and the others $0$. As in the
proof of \cite[4.3.5]{HP} we can assume that
$(U(s,N))_{s\in S}\sqcup(U(s',N))_{s'\in S'}$ forms an independent family of standard
unitary matrices thanks to the independence and assumption (2). We fix $s \in S$ and
assume $\lim_{N\rightarrow\infty} n_0(s,N)/N > 0$. (When $n(s,N)/N \to 0$ the
discussion below becomes rather trivial.) Write $A(s,N)$ and $B(s,N)$ for the matrices
appearing inside $\mathrm{Ad}\,U(s,N)$ in the above representation of $P(s,N)$,
$Q(s,N)$, that is, $A(s,N) = U(s,N)^* P(s,N) U(s,N)$ and
$B(s,N) = U(s,N)^* Q(s,N) U(s,N)$. By assumption (3) one observes that the empirical
distribution $n(s,N)^{-1}\sum_{i=1}^{n(s,N)}\delta_{x_i(s,N)}$ converges to a measure
$\rho_s$ on $[0,1]$ weakly in the almost sure sense as $N\rightarrow\infty$. Choose
(non-random) $0 \leq \xi_1(s,N) \leq \cdots \leq \xi_{n(s,N)}(s,N) \leq 1$ in such a
way that $n(s,N)^{-1}\sum_{i=1}^{n(s,N)}\delta_{\xi_i(s,N)}$ converges to $\rho_s$
weakly as $N\rightarrow\infty$. Let $\Xi(s,N)$ be the diagonal matrix with diagonal
entries $\xi_1(s,N),\dots,\xi_{n(s,N)}(s,N)$ and define 
\begin{equation*} 
C(s,N) := \begin{bmatrix} \Xi(s,N) & \sqrt{\Xi(s,N)(I-\Xi(s,N))} \\
\sqrt{\Xi(s,N)(I-\Xi(s,N))} & I-\Xi(s,N) \end{bmatrix}
\oplus I \oplus 0 \oplus I \oplus 0. 
\end{equation*} 
By \cite[4.3.4]{HP} we then have 
\begin{equation*} 
\lim_{N\rightarrow\infty}\Vert X(s,N)
- \Xi(s,N)\Vert_{p,n(s,N)^{-1}\mathrm{Tr}_{n(s,N)}} = 0
\ \ \text{almost surely\ \ for all $p\geq1$}
\end{equation*} 
so that for any polynomial $F$
\begin{equation*} 
\lim_{N\rightarrow\infty}\Vert F(B(s,N))
- F(C(s,N))\Vert_{p,N^{-1}\mathrm{Tr}_{N}} = 0
\ \ \text{almost surely\ \ for all $p\geq1$}.
\end{equation*}
Moreover, note that $(A(s,N),C(s,N))_{s\in S}$ has the limit distribution. Under these
preparations the proof is completed by the same argument as in \cite[4.3.5]{HP}.
\end{proof}

\subsection{Free transportation cost inequality for projections}
Let $\mathcal{A}^{(2n+n')}_{\mathrm{proj}}$ be the universal free product
$C^*$-algebra of $2n+n'$ copies of $C^*(\mathbb{Z}_2) = \bC\oplus\bC$,
and denote the canonical $2n+n'$ generators of projections by
$e_1,f_1,\dots,e_n,f_n,e'_1,\dots,e'_{n'}$. For a given $2n+n'$-tuple
$\vec P=(P_1,Q_1,\dots,P_n,Q_n,R_1,\dots,R_{n'})$ of projections in
$M_N(\bC)$, there is a unique $*$-homomorphism from
$\mathcal{A}^{(2n+n')}_{\mathrm{proj}}$ into $M_N(\bC)$ sending $e_i,f_i,e'_j$
to $P_i,Q_i,R_j$, respectively, which we denote by
$h \in \mathcal{A}^{(2n+n')}_{\mathrm{proj}} \mapsto h(\vec P)\in M_N(\bC)$.
For $\vec k=(k_1,l_1,\dots,k_n,l_n,k'_1,\dots,k'_{n'})\in\{0,1,\dots,N\}^{2n+n'}$,
denote by $G(N,\vec k)$ the product   
$\prod_{i=1}^n\bigl(G(N,k_i)\times G(N,l_i)\bigr)\times\prod_{j=1}^{n'}G(N,k'_j)$
of Grassmannian manifolds, and by $\mathcal{P}\bigl(G(N,\vec k)\bigr)$ the set of
Borel probability measures on $G(N,\vec k)$. Note that each
$\lambda \in \mathcal{P}\bigl(G(N,\vec k)\bigr)$ clearly gives rise to the unique
tracial state $\hat{\lambda}$ on $\mathcal{A}^{(2n+n')}_{\mathrm{proj}}$ defined by 
\begin{equation*} 
\hat{\lambda}(h) := 
\int \frac{1}{N}\mathrm{Tr}_N\bigl(h(\vec P)\bigr)\,d\lambda(\vec P)
\quad\mbox{for $h \in \mathcal{A}^{(2n+n')}_{\mathrm{proj}}$}.
\end{equation*}

Let us denote by $TS\bigl(\mathcal{A}^{(2n+n')}_{\mathrm{proj}}\bigr)$ the set of
tracial states on $\mathcal{A}^{(2n+n')}_{\mathrm{proj}}$, and moreover, for each
$\vec{\alpha} := (\alpha_1,\beta_1,\dots,\alpha_n,\beta_n,\alpha'_1,\dots,
\alpha'_{n'})\in[0,1]^{2n+n'}$, by $TS_{\vec{\alpha}}(\mathcal{A}^{(2n+n')}_{\proj})$
the set of $\tau \in TS\bigl(\mathcal{A}^{(2n+n')}_{\mathrm{proj}}\bigr)$ such that
$\tau(e_i) = \alpha_i$, $\tau(f_i) = \beta_i$ and $\tau(e'_j) = \alpha'_j$. For
$\tau_1, \tau_2 \in TS\bigl(\mathcal{A}^{(2n+n')}_{\mathrm{proj}}\bigr)$, the (free
probabilistic) Wasserstein distance $W_{2,\mathrm{free}}(\tau_1,\tau_2)$ is defined
to be the infimum of 
$$
\sqrt{\tau\Biggl(\sum_{i=1}^n\bigl(
|\sigma_1(e_i) - \sigma_2(e_i)|^2+|\sigma_1(f_i) - \sigma_2(f_i)|^2\bigr)
+\sum_{j=1}^{n'} |\sigma_1(e'_j) - \sigma_2(e'_j)|^2\Biggr)}
$$ 
over all $\tau \in TS\bigl(\mathcal{A}^{(2n+n')}_{\mathrm{proj}}
\bigstar\mathcal{A}^{(2n+n')}_{\mathrm{proj}}\bigr)$ with $\tau\circ\sigma_1=\tau_1$,
$\tau\circ\sigma_2=\tau_2$, where $\sigma_1$ and $\sigma_2$ stand for the canonical
embedding maps of $\mathcal{A}^{(2n+n')}_{\mathrm{proj}}$ into the left and right
copies in $\mathcal{A}^{(2n+n')}_{\mathrm{proj}}
\bigstar\mathcal{A}^{(2n+n')}_{\mathrm{proj}}$, respectively. The next lemma will be
one of the keys in proving a free transportation cost inequality.  

\begin{lemma}\label{L-5.2}
For each pair $\lambda_1, \lambda_2 \in \mathcal{P}\bigl(G(N,\vec k)\bigr)$ we have
\begin{equation*} 
W_{2,\mathrm{free}}(\hat{\lambda}_1,\hat{\lambda}_2)
\leq \frac{1}{\sqrt{N}} W_{2,HS}(\lambda_1,\lambda_2)
\leq \frac{1}{\sqrt{N}}W_{2,d}(\lambda_1,\lambda_2).
\end{equation*}
Here, $W_{2,HS}$ and $W_{2,d}$ are the usual Wasserstein distances determined by
the Hilbert-Schmidt norm $\Vert P-Q\Vert_{HS}$ and the geodesic distance $d(P,Q)$ with
respect to the Riemannian metric induced from $\mathrm{Tr}_N$, respectively.
\end{lemma}

\begin{proof}
The first inequality is shown in the same way as in \cite[Lemma 1.3]{HU1}, while the
second immediately follows from the inequality $\Vert P-Q\Vert_{HS} \leq d(P,Q)$ (see
e.g.~\cite[Appendix B]{He}).
\end{proof} 

Let $\vec{\alpha} \in[0,1]^{2n+n'}$ and $\vec k(N)=(k_1(N),l_1(N),\dots,k_n(N),l_n(N),
k'_1(N),\dots,k'_{n'}(N))\in\{0,1,\dots,N\}^{2n+n'}$ for $N\in\bN$ be given so that
$\vec k(N)/N\to\vec\alpha$ as $N\to\infty$. For each
$\tau \in TS_{\vec{\alpha}}\bigl(\mathcal{A}^{(2n+n')}_{\mathrm{proj}}\bigr)$ the free
entropy $\chi_{\mathrm{proj}}(\tau)$ is defined as follows. We denote by  
$\Gamma_\proj(\tau;\vec k(N);\allowbreak N,m,\varepsilon)$ the set of all
$2n+n'$-tuples $\vec P\in G(N,\vec k(N))$ such that
\begin{equation*} 
\left|\frac{1}{N}\mathrm{Tr}_N\bigl(h(\vec P)\bigr) - \tau(h)\right| < \varepsilon
\end{equation*}
for all monomials $h \in \mathcal{A}^{(2n+n')}_{\proj}$ in
$e_1,f_1,\dots,e_n,f_n,r_1,\dots,r_{n'}$ of degree at most $m$. We then define 
\begin{equation*} 
\chi_{\mathrm{proj}}(\tau)
:= \lim_{\substack{m\rightarrow\infty \\ \varepsilon \searrow 0}}
\limsup_{N\rightarrow\infty}\,\frac{1}{N^2}\log
\gamma_{\vec k(N)}\bigl(\Gamma_\proj(\tau;\vec k(N);N,m,\varepsilon)\bigr), 
\end{equation*}
where $\gamma_{\vec k(N)}:=\bigotimes_{i=1}^n\bigl(\gamma_{G(N,k_i(N))}\otimes
\gamma_{G(N,l_i(N))}\bigr)\otimes\bigotimes_{j=1}^{n'}\gamma_{G(N,k'_j(N))}$ on
$G(N,\vec k(N))$. Note that the quantity $\chi_{\proj}(\tau)$ is noting less than
$\chi_{\proj}(p_1,q_1,\dots,p_n,q_n,r_1,\dots,r_{n'})$ when $p_i:=\pi_\tau(e_i)$,
$q_i:=\pi_\tau(f_i)$ and $r_j:=\pi_\tau(e'_j)$ in the GNS representation of
$\mathcal{A}_{\proj}^{(2n+n')}$ associated with $\tau$; hence it is independent of
the particular choice of $\vec k(N)$ due to Proposition \ref{P-1.1}.  

In what follows, let
$\tau \in TS_{\vec{\alpha}}\bigl(\mathcal{A}^{(2n+n')}_{\mathrm{proj}}\bigr)$ be
arbitrarily fixed. Then one can choose a subsequence $N_1 < N_2 < \cdots$ so that 
\begin{equation}\label{F-5.1}
\chi_{\proj}(\tau) = \lim_{m\rightarrow\infty}\frac{1}{N_m^2}\log
\gamma_{\vec k(N_m)}\bigl(\Gamma_\proj(\tau;\vec k(N_m);N_m,m,1/m)\bigr).
\end{equation}
Set $\Gamma_{N_m}:=\Gamma_\proj(\tau;\vec k(N_m);N_m,m,1/m)$ and define
$\lambda_{N_m}^{\tau} \in \mathcal{P}\bigl(G(N_m,\vec k(N_m))\bigr)$ by 
\begin{equation*} 
d\lambda_{N_m}^{\tau}(\vec P) := \frac{1}{\gamma_{\vec k(N_m)}(\Gamma_{N_m})}
\,\mathbf{1}_{\Gamma_{N_m}}(\vec P)\,d\gamma_{\vec k(N_m)}(\vec P). 
\end{equation*}  

\begin{lemma}\label{L-5.3}
$\displaystyle \lim_{m\rightarrow\infty} \hat{\lambda}^{\tau}_{N_m} = \tau$ in the
weak* topology.
\end{lemma}

\begin{proof} The proof can be found in \cite{HU1}, even though only the self-adjoint
and the unitary cases are treated there.
\end{proof} 

For $1\le i\le n$ the $C^*$-subalgebra generated by $e_i,f_i$ (obviously identified
with $\cA^{(2)}_\proj=C^*(\mathbb{Z}_2\bigstar\mathbb{Z}_2)$) is isomorphic to 
\begin{equation*}  
\mathcal{A} := \bigl\{ a(\cdot) = [a_{ij}(\cdot)]_{i,j=1}^2 \in C([0,1];M_2(\bC)) :
\text{$a(0)$, $a(1)$ are diagonals} \bigr\}
\end{equation*} 
by the $*$-isomorphism given by
\begin{equation*} 
e_i \mapsto e(t) = \begin{bmatrix} 1 & 0 \\ 0 & 0 \end{bmatrix}, \quad 
f_i \mapsto f(t) = \begin{bmatrix} t & \sqrt{t(1-t)} \\ \sqrt{t(1-t)} & 1-t \end{bmatrix},   
\end{equation*}  
and any tracial state on $\cA$ is of the form
\begin{equation*} 
\tau_{\nu,\{\alpha_{ij}\}}(a) :=  
\alpha_{10} a_{11}(0) + \alpha_{01}a_{22}(0) + \alpha_{11} a_{11}(1)
+ \alpha_{00} a_{22}(1) + \int_{(0,1)} \frac{1}{2}\mathrm{Tr}_2(a(t))\,d\nu(t),  
\end{equation*} 
where $\alpha_{ij}\geq 0$, $\sum_{i,j=0}^1\alpha_{ij} \leq 1$ and $\nu$ is a measure
on $(0,1)$  with $\nu((0,1)) = 1-\sum_{i,j=0}^1\alpha_{ij}$. 
Let $\vec{\psi} = (\psi_1,\dots,\psi_n)$ be an $n$-tuple of continuous functions on
$[0,1]$, and for $1\le i\le n$ define the probability distribution
$\lambda_N^{\psi_i}$ on $G(N,k_i(N))\times G(N,l_i(N))$ by 
\begin{equation*} 
d\lambda_N^{\psi_i}(P,Q) := \frac{1}{Z^{\psi_i}_N}
\exp\bigl(-N\mathrm{Tr}_N(\psi_i(PQP))\bigr)
\,d(\gamma_{G(N,k_i(N))}\otimes\gamma_{G(N,l_i(N))})(P,Q) 
\end{equation*} 
with the normalization constant $Z^{\psi_i}_N$. For $\tau_{\nu,\{\alpha_{ij}\}} \in
TS_{(\alpha_i,\beta_i)}\bigl(\mathcal{A}_{\proj}^{(2)})$ one has  
\begin{align*} 
&\chi_{\proj}(\tau_{\nu,\{\alpha_{ij}\}})
- \tau_{\nu,\{\alpha_{ij}\}}(\psi_i(efe)) \\
&\qquad= 
\frac{1}{4}\Sigma(\nu) + \frac{1}{2}\int_{(0,1)}\bigl((\alpha_{01}+\alpha_{10})\log t
+ (\alpha_{00}+\alpha_{11})\log(1-t) -\psi_i(t)\bigr)\,d\nu(t) - C    
\end{align*}   
with some constant $C$ if $\alpha_{00}\alpha_{11} = \alpha_{01}\alpha_{10}=0$, and
otherwise $-\infty$. Thus, a general result on weighted logarithmic energy (see 
\cite{ST}) ensures that there is a unique maximizer
$\tau_{(\alpha_i,\beta_i)}^{\psi_i} \in
TS_{(\alpha_i,\beta_i)}\bigl(\mathcal{A}_{\proj}^{(2)}\bigr)$ of the functional
$\tau \in TS_{(\alpha_i,\beta_i)}\bigl(\mathcal{A}_{\proj}^{(2)}\bigr) \mapsto
\chi_{\proj}(\tau) - \tau(\psi_i(efe))$. Then, we define the tracial state
$\tau_{\vec{\alpha}}^{\vec\psi} \in TS\bigl(\mathcal{A}_{\proj}^{(2n+n'}\bigr)$ by 
\begin{equation*} 
\tau_{\vec{\alpha}}^{\vec\psi}
:= \Bigl(\bigstar_{i=1}^n \tau_{(\alpha_i,\beta_i)}^{\psi_i}\Bigr) \bigstar
\tau_{(\alpha'_1,\dots,\alpha'_{n'})},
\quad\tau_{(\alpha'_1,\dots,\alpha'_{n'})}
:=\bigstar_{j=1}^{n'}\bigl(\alpha'_j\delta_0 + (1-\alpha'_j)\delta_1\bigr)
\end{equation*} 
in the natural identification $\mathcal{A}_{\proj}^{(2n+n')}
= \Bigl(\bigstar_{i=1}^n \mathcal{A}_{\proj}^{(2)}\Bigr)\bigstar
\Bigl(\bigstar_{j=1}^{n'}C^*(\bZ_2)\Bigr)$. Furthermore, we define the joint distribution
\begin{equation*}
\lambda_N^{\vec\psi} := \Biggl(\bigotimes_{i=1}^n \lambda_N^{\psi_i}\Biggr)\otimes
\Biggl(\bigotimes_{j=1}^{n'} \gamma_{G(N,k'_j(N))}\Biggr)
\quad\mbox{on $G(N,\vec k(N))$}
\end{equation*} 
(also considered as a $2n+n'$-tuple of random projection matrices). The next lemma
follows from a large deviation result for two projection matrices in \cite{HP1} and
Theorem \ref{T-5.1}.     

\begin{lemma}\label{L-5.4} \ 
\begin{enumerate} 
\item $\displaystyle B_{\psi_i} := \lim_{N\rightarrow\infty}\frac{1}{N^2}
\log Z_N^{\psi_i}$ exists for every $1\le i\le n$.   
\item $\displaystyle \lim_{N\rightarrow\infty}\hat{\lambda}_N^{\vec\psi}
= \tau_{\vec{\alpha}}^{\vec\psi}$ in the weak* topology.
\end{enumerate}  
\end{lemma}

\begin{proof}
When $p_1,q_1,\dots,p_n,q_n$ are absent, we have nothing to do for (1) and moreover
(2) immediately follows from Voiculescu's original result \cite[Theorem 3.11]{V0}
rather than Theorem \ref{T-5.1} as follows. Let $R_1(N),\dots,R_{n'}(N)$ be an
independent family of random projection matrices of ranks $k'_j(N)$ distributed under
$\gamma_{G(N,k'_j(N))}$, respectively. Note that $\lambda_N^{\vec\psi}$ in this
case is nothing less than
$\gamma_N := \bigotimes_{j=1}^{n'}\gamma_{G(N,k'_j(N))}$. For a monomial
$h = r_{j_1}\cdots r_{j_m} \in \mathcal{A}_{\proj}^{(n')}$ one has 
\begin{align*} 
\hat\gamma_N(h) = \mathbb{E}\circ\biggl(\frac{1}{N}\mathrm{Tr}_N\biggr)
(R_{j_1}(N)\cdots R_{j_m}(N)), 
\end{align*}
which converges to $\tau_{(\alpha'_1,\dots,\alpha'_{n'})}(r_{j_1}\cdots r_{j_m})$
thanks to \cite[Theorem 3.11]{V0}. This immediately implies (2) in this special case. 

For the general case, i.e., when $p_1,q_1,\dots,p_n,q_n$ really appear, we need to
show (1) and $\lim_{N\rightarrow\infty} \hat{\lambda}_N^{\psi_i}
= \tau_{\vec{\alpha}}^{\psi_i}$ weakly* for each $1\le i\le n$. Both are simple
applications of the large deviation result for the empirical eigenvalue distribution
of two random projection matrix pair $(P(N),Q(N))$ distributed under
$\lambda_N^{\psi_i}$, whose proof is essentially same as in \cite[Proposition 2.1]{HP1} 
(or in the proof of \cite[Theorem 2.1]{HU2}). Once the latter convergence was
established, the above argument would equally work well even in the general setting
when \cite[Theorem 3.11]{V0} is replaced by Theorem \ref{T-5.1}.    
\end{proof} 

With the above lemmas we can now prove the following transportation cost inequality 
in the essentially same manner as in \cite{HU1}.  

\begin{thm}\label{T-5.4} Assume that $\psi_i$'s are $C^2$-functions and $\rho :=
\min\bigl\{1 - c_1\Vert \psi'_i\Vert_{\infty} - c_2\Vert\psi''_i\Vert_{\infty}:
1\leq i\leq n\bigr\} > 0$  for some universal constants $c_1, c_2 > 0$ {\rm (}for
example, one can choose $c_1 = 6$, $c_2 = 9/2$, but these do not seem optimal\,{\rm )}.
For every $\tau \in TS_{\vec{\alpha}}\bigl(\mathcal{A}_{\proj}^{(2n+n')}\bigr)$ we have 
\begin{equation}\label{F-5.2} 
W_{2,\mathrm{free}}\Bigl(\tau,\tau_{\vec{\alpha}}^{\vec\psi}\Bigr) \leq 
\sqrt{\frac{2}{\rho}\Biggl(-\chi_{\proj}(\tau)
+ \tau\Biggl(\sum_{i=1}^n \psi_i(pqp)\Biggr) + B_{\vec\psi}\Biggr)}
\end{equation} 
with $B_{\vec\psi} := \sum_{i=1}^n B_{\psi_i}$. In particular, when
$p_1,q_1,\dots,p_n,q_n$ are absent, \eqref{F-5.2} simply becomes 
\begin{equation}\label{F-5.3} 
W_{2,\mathrm{free}}\bigl(\tau,\tau_{(\alpha'_1,\dots,\alpha'_{n'})}\bigr)
\leq \sqrt{-2\chi_{\mathrm{proj}}(\tau)}. 
\end{equation} 
\end{thm}

\begin{proof} 
Since $W_{2,\mathrm{free}}$ is lower semi-continuous in the weak* topology, we have
by Lemmas \ref{L-5.3} and \ref{L-5.4}\,(2) 
\begin{equation*}
W_{2,\mathrm{free}}\Bigl(\tau,\tau_{\vec{\alpha}}^{\vec\psi}\Bigr)
\leq \liminf_{m\rightarrow\infty} 
W_{2,\mathrm{free}}\Bigl(\hat{\lambda}^{\tau}_{N_m},\hat{\lambda}^{\vec\psi}_{N_m}\Bigr). 
\end{equation*} 
By Lemma \ref{L-5.2} we also have 
\begin{equation*}
W_{2,\mathrm{free}}\Bigl(\hat{\lambda}^{\tau}_{N_m},
\hat{\lambda}^{\vec\psi}_{N_m}\Bigr) 
\leq \frac{1}{\sqrt{N_m}} 
W_{2,d}\Bigl(\lambda^{\tau}_{N_m},\lambda^{\vec\psi}_{N_m}\Bigr). 
\end{equation*} 
We then need to confirm Bakry and Emery's $\Gamma_2$-criterion \cite{BE} for
$\lambda^{\vec\psi}_N$ with the constant $\rho N$, that is, 
\begin{equation}\label{F-5.4} 
\mathrm{Ric}\bigl(G(N,\vec k(N))\bigr) + \mathrm{Hess}(\Psi_N) \geq \rho NI_{d(N)},
\end{equation} 
where $\mathrm{Ric}\bigl(G(N,\vec k(N))\bigr)$ is the Ricci curvature tensor of
$G(N,\vec k(N))$, $\mathrm{Hess}(\Psi_N)$ is the Hessian of the trace function 
\begin{align*} 
\Psi_N(P_1,Q_1,\dots,P_n,Q_n,R_1,\dots,R_{n'}) := 
N\mathrm{Tr}_N\Biggl(\sum_{i=1}^n \psi_i(P_i Q_i P_i)\Biggr), 
\end{align*}
and $d(N)$ is the dimension of $G(N,\vec k(N))$, i.e.,
$$
d(N):=2\sum_{i=1}^n\bigl(k_i(N)(N-k_i(N))+l_i(N)(N-l_i(N))\bigr)
+2\sum_{j=1}^{n'}k'_j(N)(N-k'_j(N)).
$$
It is known (see \cite[Eq.\,(2.2)]{HU2}) that $\mathrm{Ric}(G(N,k)) = NI_{2k(N-k)}$
so that we need only to estimate the Hessian $\mathrm{Hess}(\Psi_N^{(i)})$ of the
trace function $\Psi_N^{(i)} : (P,Q) \in G(N,k_i(N))\times G(N,l_i(N)) \mapsto
N\mathrm{Tr}_N(\psi_i(PQP))$ from below for each $1\le i\le n$. This can be done by
computing $\mathrm{Hess}(\Psi_N^{(i)})$ explicitly, 
and consequently we can find two universal constants $c_1, c_2 > 0$ so that
\begin{equation*}
\mathrm{Hess}(\Psi_N^{(i)})
\geq -N(c_1\Vert \psi'_i\Vert_{\infty} + c_2\Vert \psi''_i\Vert_{\infty})
I_{2k_i(N)(N-k_i(N))+2l_i(N)(N-l_i(N))}.
\end{equation*}
Hence \eqref{F-5.4} is confirmed. See Remark \ref{R-5.6} below for more details on
this estimate. Thus, by the transportation cost inequality in the Riemannian manifold
setting due to Otto and Villani \cite{OV} we obtain
\begin{equation}\label{F-5.5}
W_{2,d}\Bigl(\lambda^{\tau}_{N_m},\lambda^{\vec\psi}_{N_m}\Bigr)
\leq \sqrt{\frac{2}{\rho N_m}
S\Bigl(\lambda^{\tau}_{N_m},\lambda^{\vec\psi}_{N_m}\Bigr)},
\end{equation}
where $S\Bigl(\lambda^{\tau}_{N_m},\lambda^{\vec\psi}_{N_m}\Bigr)$ stands for the
usual relative entropy. We compute  
\begin{align*}
S\Bigl(\lambda^{\tau}_{N_m},\lambda^{\vec\psi}_{N_m}\Bigr) 
&=\int\log\frac{d\lambda^{\tau}_{N_m}}{d\lambda^{\vec\psi}_{N_m}}
\,\lambda^{\tau}_{N_m} \\
&=-\log\gamma_{\vec k(N_m)}(\Gamma_{N_m})
+ N_m^2 \hat{\lambda}_{N_m}^{\tau}\Biggl(\sum_{i=1}^n \psi_i(e_i f_i e_i)\Biggr)
+ \sum_{i=1}^n\log Z_{N_m}^{\psi_i}.  
\end{align*}   
Consequently, we obtain the desired inequality \eqref{F-5.2} by taking the limit of
\eqref{F-5.5} as $m\rightarrow\infty$ after divided by $N_m^2$ due to \eqref{F-5.1},
Lemmas \ref{L-5.3} and \ref{L-5.4}\,(1). Finally, we should remark that if
$p_1,q_1,\dots,p_n,q_n$ disappeared, then the argument would become simpler without
estimating the Hessian of $\Psi_N$.    
\end{proof}  

\begin{remark}\label{R-5.6}{\rm
The computation of $\mathrm{Hess}(\Psi_N^{(i)})$ mentioned in the above
proof is outlined here. The tangent space $T_P G(N,k)$ at $P \in G(N,k)$ is identified
with the set of $X \in M_N(\mathbb{C})^{sa}$ satisfying $X = PX + XP$, on which our
Riemannian metric is given by $\langle X|Y \rangle := \mathrm{Re}\,\mathrm{Tr}_N(YX)$
(this is inherited from that on the Euclidean space $M_N(\mathbb{C})^{sa}$). Moreover,
the geodesic curve started at $P$ with tangent vector $X$ is given by
$C(t) := \exp(t[X,P])P\exp(-t[X,P])$ for $t \in \mathbb{R}$. (See e.g.~\cite[\S2]{CPR}
for a brief summary and references therein.) Since 
\begin{align*} 
&\bigl\langle\mathrm{Hess}(\Psi^{(i)}_N)((C_1(0),C_2(0))(C_1'(0)\oplus C_2'(0)
|C_1'(0)\oplus C_2'(0))\bigr\rangle \\
&\qquad= \frac{d^2}{dt^2}\bigg|_{t=0}N\mathrm{Tr}_N(\psi_i(C_1(t)C_2(t)C_1(t)))
\end{align*} 
for geodesic curves $C_1(t) \in G(N,k_i(N))$ and $C_2(t) \in G(N,l_i(N))$, it suffices
{\rm (}for getting the desired inequality in the above proof{\rm )} to estimate, at
$t=0$, the second derivative of the composition of
$\phi(t) := C_1(t)C_2(t)C_1(t) \in M_N(\mathbb{C})^{sa}$ and
$X \in M_N(\mathbb{C})^{sa} \mapsto \Phi(X) := N\mathrm{Tr}_N(\psi_i(X))$ with the
usual Euclidean structure on $M_N(\mathbb{C})^{sa}$. Passing once to the identification
$M_N(\mathbb{C})^{sa} = \mathbb{R}^{N^2}$, we observe that  
\begin{equation*}
(\Phi\circ\phi)''(0)
= \bigl\langle(\nabla^2 \Phi)(\phi(0))\phi'(0)|\phi'(0)\bigr\rangle
+ \bigl\langle(\nabla \Phi)(\phi(0))|\phi''(0)\bigr\rangle 
\end{equation*}
thanks to the usual chain rule. By \cite[Lemma 1.2]{HPU} we can estimate the operator
norms $\Vert (\nabla^2\Phi)(\phi(0)) \Vert_{\infty}$ (for linear operators on
$(M_N(\mathbb{C})^{sa},\langle\,\cdot\,|\,\cdot\,\rangle)$) and  
$\Vert (\nabla\Phi)(\phi(0))\Vert_{\infty}$ (for elements in $M_N(\mathbb{C})^{sa}$)
by $N\Vert \psi_i''\Vert_{\infty}$  and $N\Vert \psi_i'\Vert_{\infty}$, respectively,
from the above. As mentioned above the tangent vector
$C_i'(0) \in M_N(\mathbb{C})^{sa}$ satisfies $C_i'(0) = C_i(0)C_i'(0) + C_i'(0)C_i(0)$
and the geodesic curve $C_i(t)$ must be 
\begin{equation*} 
C_i(t) = \exp(t[C_i'(0),C_i(0)])C_i(0)\exp(-t[C_i'(0),C_i(0)]). 
\end{equation*} 
It follows from these facts that 
\begin{align*} 
\phi'(0) 
&= 
C'_1(0)C_2(0)C_1(0) + C_1(0)C_2'(0)C_1(0) + C_1(0)C_2(0)C_1'(0), \\
\phi''(0)
&= 
[[C_1'(0),C_1(0)],C_1'(0)]C_2(0)C_1(0) \\
&\quad+ C_1(0)[[C_2'(0),C_2(0)],C_2'(0)]C_1(0) \\
&\quad\quad+ C_1(0)C_2(0)[[C_1'(0),C_1(0)],C_1'(0)] \\
&\quad\quad\quad+ 
2\big\{C_1'(0)C_2'(0)C_1(0) + C_1'(0)C_2(0)C_1'(0) + C_1(0)C_2'(0)C_1'(0)\big\}.    
\end{align*}
Hence we get the rough estimates
\begin{align*} 
\Vert\phi'(0)\Vert_{HS}^2
&\leq 6\Vert C_1'(0)\Vert_{HS}^2 + 3\Vert C_2'(0)\Vert_{HS}^2, \\
\Vert\phi''(0)\Vert_{1,\mathrm{Tr}_N}
&\leq 8\Vert C_1'(0)\Vert_{HS}^2 + 4\Vert C_2'(0)\Vert_{HS}^2 
\end{align*}  
(we used $2C_i'(0)^2 = |[C_i'(0),C_i(0)],C_i'(0)]|$, $i=1,2$, for the latter).
Therefore,
\begin{equation*} 
(\Phi\circ\phi)''(0)
\leq N\big\{(8\Vert\psi_i'\Vert_{\infty} + 6\Vert\psi_i''\Vert_{\infty})
\Vert C_1'(0)\Vert_{HS}^2
+ (4\Vert\psi_i'\Vert_{\infty} + 3\Vert\psi_i''\Vert_{\infty})
\Vert C_1'(0)\Vert_{HS}^2\big\}. 
\end{equation*} 
Since $\Phi\circ\phi(t)$ does not change when $C_1(t), C_2(t)$ are interchanged, one
finally finds two universal constants $c_1 = 6>0$, $c_2 = 9/2 > 0$ so that
$$
|(\Phi\circ\phi)''(0)| \leq N(c_1\Vert\psi_i'\Vert_{\infty}
+c_2\Vert\psi_i''\Vert_{\infty})(\Vert C_1'(0) \Vert_{HS}^2 + \Vert C_2'(0)
\Vert_{HS}^2),
$$
which immediately implies the desired inequality.

Finally, it should be pointed out that $(6,9/2)$ can be also chosen for two universal
constants $(c_1,c_2)$ in \cite[Proposition 3.1]{HU2}.}   
\end{remark} 

\begin{cor}\label{C-5.5}
If $p_1,\dots,p_n$ are projections in $(\cM,\tau)$ and
$\chi_{\mathrm{proj}}(p_1,\dots,p_n) = 0$, then $p_1,\dots,\allowbreak p_n$ are free.    
\end{cor}

\begin{proof}
This follows from \eqref{F-5.3} and the fact that $W_{2,\mathrm{free}}$ is a metric
on $TS_{\vec\alpha}\bigl(\cA^{(n)}_\proj\bigr)$ where
$\vec\alpha:=(\tau(p_1),\dots,\tau(p_n))$. 
\end{proof} 

The corollary was an essential ingredient of the proof of Theorem \ref{T-4.1}. In the
self-adjoint case, the free transportation cost inequality \cite[Theorem 2.2 or
Corollary 2.3]{HU1} provides a new proof of the fact that $X_1,\dots,X_n$ form a free
semicircular system if $\chi(X_1,\dots,X_n)$ attains the maximum under the restriction
$\tau(X_i^2)=1$, while Voiculescu's original proof in \cite{V4} is based on the
infinitesimal change of variable formula. 

\section{Free pressure}
\setcounter{equation}{0}

Let $\bigl(\cA^{(n)}_\proj\bigr)^{sa}$ denote the space of self-adjoint elements in
the universal $C^*$-algebra $\cA^{(n)}_\proj$ with the canonical projection generators
$e_1,\dots,e_n$ as given in the previous section. Elements in
$\bigl(\cA^{(n)}_\proj\bigr)^{sa}$ are considered as ``free probabilistic hamiltonians
on $\bZ_2^{\bigstar n}$." Motivated from the statistical mechanical  viewpoint, we
introduce the free pressure for those free hamiltonians, and its Legendre
transform with respect to the duality between $\bigl(\cA^{(n)}_\proj\bigr)^{sa}$ and
$TS\bigl(\cA^{(n)}_\proj\bigr)$ is compared with $\chi_\proj$.

Let $\vec\alpha=(\alpha_1,\dots,\alpha_n)\in[0,1]^n$ and
$\vec k(N)=(k_1(N),\dots,k_n(N))\in\{0,1,\dots,N\}$ for $N\in\bN$ be given so that
$\vec k(N)/N\to{\vec\alpha}$ as $N\to\infty$. As before we write
$G(N,\vec k(N)):=\prod_{i=1}^nG(N,k_i(N))$ and
$\gamma_{\vec k(N)}:=\bigotimes_{i=1}^n\gamma_{G(N,k_i(N))}$ for short. For
$\vec P=(P_1,\dots,P_n)\in G(N)$ we have the $*$-homomorphism
$h\in\cA^{(n)}_\proj\mapsto h(\vec P)\in M_N(\bC)$ sending $e_i$ to
$P_i$, $1\le i\le n$. For each $h\in\bigl(\cA^{(n)}_\proj\bigr)^{sa}$ define
\begin{equation}\label{F-6.1}
\pi_{\vec\alpha}(h):=\limsup_{N\to\infty}
{1\over N^2}\log\int_{G(N,\vec k(N))}\exp\Bigl(
-N\Tr_N\bigl(h(\vec P)\bigr)\Bigr)\,d\gamma_{\vec k(N)}(\vec P),
\end{equation}
which we call the {\it free pressure} of $h$ under the trace values
$(\alpha_1,\dots,\alpha_n)$.

\begin{prop}\label{P-6.1}
The above definition of $\pi_{\vec\alpha}(h)$ is independent of the choices of
$\vec k(N)$ with $\vec k(N)/N\to\vec\alpha$.
\end{prop}

\begin{proof}
Let $\vec l(N)=(l_1(N),\dots,l_n(N)))$, $N\in\bN$, be another sequence such that
$\vec l(N)/N\to\vec\alpha$ as $N\to\infty$. For $h\in\cA^{(n)}_\proj$ and $N\in\bN$
set
$$
\delta_N(h):=\max_{\vec U\in\U(N)^n}
\bigg|{1\over N}\Tr_N\bigl(h(\xi_{\vec k(N)}(\vec U)\bigr)
-{1\over N}\Tr_N\bigl(h(\xi_{\vec l(N)}(\vec U)\bigr)\bigg|,
$$
where $\xi_{\vec k(N)}(\vec U):=(\xi_{N,k_1(N)}(U_1),\dots,\xi_{N,k_n(N)}(U_n))$ for
$\vec U=(U_1,\dots,U_n)$ (see \S1). For each $h\in\bigl(\cA^{(n)}_\proj\bigr)^{sa}$ we
get thanks to \eqref{F-1.2}
\begin{align*}
&\Bigg|{1\over N^2}\log\int_{G(N,\vec k(N))}
\exp\Bigl(-N\Tr_N\bigl(h(\vec P)\bigr)\Bigr)
\,d\gamma_{\vec k(N)}(\vec P) \\
&\qquad\quad-{1\over N^2}\log\int_{G(N,\vec k(N))}
\exp\Bigl(-N\Tr_N\bigl(h(\vec P)\bigr)\Bigr)
\,d\gamma_{\vec k(N)}(\vec P)\Bigg| \\
&\quad=\Bigg|{1\over N^2}\log\int_{\U(N)^n}\exp\Bigl(-N\Tr_N\bigl(
h(\xi_{\vec k(N)}(\vec U)\bigr)\Bigr)
\,d\bigl(\gamma_{\U(N)}\bigr)^{\otimes n}(\vec U) \\
&\qquad\quad-{1\over N^2}\log\int_{\U(N)^n}\exp\Bigl(-N\Tr_N\bigl(
h(\xi_{\vec l(N)}(\vec U)\bigr)\Bigr)
\,d\bigl(\gamma_{\U(N)}\bigr)^{\otimes n}(\vec U)\Bigg| \\
&\quad\le\delta_N(h).
\end{align*}
Hence, it suffices to prove that $\delta_N(h)\to0$ as $N\to\infty$ for any
$h\in\cA^{(n)}_\proj$. Since $|\delta_N(h_1)-\delta_N(h_2)|\le\|h_1-h_2\|$ for all
$h_1,h_2\in\cA^{(n)}_\proj$, we may show that $\delta_N(h)\to0$ for $h=e_{i_1}\cdots e_{i_r}$ with $1\le
i_1,\dots,i_r\le n$. For such $h$, as in the proof of Proposition \ref{P-1.1} we have
\begin{align*}
\bigg|{1\over N}\Tr_N\bigl(h(\xi_{\vec k(N)}(\vec U)\bigr)
-{1\over N}\Tr_N\bigl(h(\xi_{\vec l(N)}(\vec U)\bigr)\bigg|
&\le\sum_{j=1}^r\|\xi_{N,k_{i_j}(N)}(U_{i_j})-\xi_{N,l_{i_j}(N)}(U_{i_j})\|_1 \\
&\le\sum_{j=1}^r{|k_{i_j}(N)-l_{i_j}(N)|\over N}\longrightarrow0
\end{align*}
as $N\to\infty$, and the conclusion follows.
\end{proof}

The following are basic properties of $\pi_{\vec\alpha}(h)$; we omit the proofs very
similar to those of \cite[Proposition 2.3]{Hi} but note that the last assertion of (iv)
follows from \eqref{F-6.2} and Proposition \ref{P-6.4}\,(1) below.

\begin{prop}\label{P-6.2} \
\begin{itemize}
\item[(i)] $\pi_{\vec\alpha}(h)$ is convex on $\bigl(\cA^{(n)}_\proj\bigr)^{sa}$.
\item[(ii)] If $h_1,h_2\in\bigl(\cA^{(n)}_\proj\bigr)^{sa}$ and $h_1\le h_2$, then
$\pi_{\vec\alpha}(h_1)\ge\pi_{\vec\alpha}(h_2)$.
\item[(iii)] $|\pi_{\vec\alpha}(h_1)-\pi_{\vec\alpha}(h_2)|\le\|h_1-h_2\|$ for all
$h_1,h_2\in\bigl(\cA^{(n)}_\proj\bigr)^{sa}$.
\item[(iv)] If $h_1\in\bigl(\cA^{(j)}_\proj\bigr)^{sa}$ and
$h_2\in\bigl(\cA^{(n-j)}_\proj\bigr)^{sa}$ with $1\le j<n$, then
$$
\pi_{\vec\alpha}(h_1+h_2)\le\pi_{(\alpha_1,\dots,\alpha_j)}(h_1)
+\pi_{(\alpha_{j+1},\dots,\alpha_n)}(h_2)
$$
where $h_1+h_2$ is the sum as an element of
$\cA^{(n)}_\proj=\cA^{(j)}_\proj\bigstar\cA^{(n-j)}_\proj$. In particular when
$j=1$ or $j=2$, equality holds in the above inequality.
\end{itemize}
\end{prop}

\begin{remark}\label{R-6.3}{\rm
Another possible definition of free pressure is to use the probability measure
$\gamma_{G(N)}^{(1)}$ or $\gamma_{G(N)}^{(2)}$ on $G(N)$ given in Remark \ref{R-1.3}.
For $h\in\bigl(\cA^{(n)}_\proj\bigr)^{sa}$ define
$$
\pi^{(j)}(h):=\limsup_{N\to\infty}{1\over N^2}\log
\int_{G(N)^n}\exp\Bigl(-N\Tr_N\bigl(h(\vec P)\bigr)\Bigr)
\,d\bigl(\gamma_{G(N)}^{(j)}\bigr)^{\otimes n}(\vec P)
$$
for $j=1,2$. It is not difficult to show that
$$
\pi^{(1)}(h)=\pi^{(2)}(h)
=\max\bigl\{\pi_{\vec\alpha}(h):\vec\alpha\in[0,1]^n\bigr\}
$$
for every $h\in\bigl(\cA^{(n)}_\proj\bigr)^{sa}$. We simply write $\pi(h)$ for these
equal quantities; then $\pi(h)$ has the same properties as in Proposition \ref{P-6.2}.
However, unlike the free entropy quantities $\chi^{(j)}_\proj$ discussed
in Remark \ref{R-1.3}, $\pi(h)$ does not coincide with $\pi_{\vec\alpha}(h)$; the
latter actually depends on $\vec\alpha$.
}\end{remark}

In the single projection case, $\cA^{(1)}_\proj=\bC^2$,
$\bigl(\cA^{(1)}_\proj\bigr)^{sa}=\bR^2$ and
$TS\bigl(\cA^{(1)}_\proj\bigr)=\{\tau_\alpha:0\le\alpha\le1\}$ where
$\tau_\alpha(\zeta_1,\zeta_2)=\alpha\zeta_1+(1-\alpha)\zeta_2$ for
$(\zeta_1,\zeta_2)\in\bC^2$. Let $0\le\alpha\le1$ and choose $k(N)$ such that
$k(N)/N\to\alpha$. For each $h=(h_1,h_2)\in\bR^2$, it is straightforward to check that
\begin{align}\label{F-6.2}
\pi_\alpha(h)&=\lim_{N\to\infty}{1\over N^2}\log
\int_{G(N,k(N))}\exp\bigl(-N\Tr_N(h(P))\bigr)\,d\gamma_{G(N,k(N))}(P) \nonumber\\
&=-\alpha h_1-(1-\alpha)h_2=-\tau_\alpha(h)
\end{align}
and hence $\chi_\proj(\tau_\alpha)=0=\tau_\alpha(h)+\pi_\alpha(h)$. Moreover,
$$
\pi(h)=-\min\{h_1,h_2\}=\min\{-\tau_\alpha(h)+\chi(\tau_\alpha):0\le\alpha\le1\}.
$$

In the case of two projections, $\cA^{(2)}_\proj=C^*(\bZ_2\bigstar\bZ_2)$ with the
canonical projection generators $e,f$. Let $\alpha,\beta\in[0,1]$. The next theorem
says that the free entropy
$\chi_\proj(\tau)$ for $\tau\in TS_{(\alpha,\beta)}\bigl(\cA^{(2)}_\proj\bigr)$ and
the free pressure $\pi_{(\alpha,\beta)}(h)$ for $h\in\bigl(\cA^{(2)}_\proj\bigr)^{sa}$
are the Legendre transforms of each other.

\begin{prop}\label{P-6.4} \
\begin{itemize}
\item[(1)] In the definition of $\pi_{(\alpha,\beta)}(h)$ in \eqref{F-6.1} $\limsup$
can be replaced by $\lim$.
\item[(2)] $\pi_{(\alpha,\beta)}(h)=\max\Bigl\{-\tau(h)+\chi_\proj(\tau):
\tau\in TS_{(\alpha,\beta)}\bigl(\cA^{(2)}_\proj\bigr)\Bigr\}$ for every
$h\in\bigl(\cA^{(2)}_\proj\bigr)^{sa}$.
\item[(3)] $\chi_\proj(\tau)=\inf\Bigl\{\tau(h)+\pi_{(\alpha,\beta)}(h):
h\in\bigl(\cA^{(2)}_\proj\bigr)^{sa}\Bigr\}$ for every
$\tau\in TS_{(\alpha,\beta)}\bigl(\cA^{(2)}_\proj\bigr)$.
\item[(4)] $\pi(h)=\max\Bigl\{-\tau(h)+\chi_\proj(\tau):
\tau\in TS\bigl(\cA^{(2)}_\proj\bigr)\Bigr\}$ for every
$h\in\bigl(\cA^{(2)}_\proj\bigr)^{sa}$.
\end{itemize}
\end{prop}

\begin{proof}
Thanks to the Lipschitz continuity in $h$ of the quantity inside $\limsup$ in
\eqref{F-6.1} as well as both sides of the equality in (2), to prove (1) and (2), we
may assume that $h$ is a self-adjoint polynomial of $e,f$ written as
$$
h=C\1+Ae+Bf+\sum_{k=1}^mA_j(efe)^j+\sum_{j=1}^mB_j(fef)^j
+\sum_{j=1}^mD_j((ef)^j+(fe)^j)
$$
with $A,B,C,A_j,B_j,D_j\in\bR$. Set
$$
h_0:=Ae+Bf+\sum_{k=0}^mC_k(efe)^k
$$
with $C_0:=C$, $C_j:=A_j+B_j+D_j$, $1\le j\le m$. We then get $\tau(h)=\tau(h_0)$ and
$\Tr_N(h(P,Q))=\Tr_N(h_0(P,Q))$ for $P,Q\in G(N)$ so that
$\pi_{(\alpha,\beta)}(h)=\pi_{(\alpha,\beta)}(h_0)$. Hence it is enough to prove (1)
and (2) for $h_0$ above. A bit more generally, let
$h\in\bigl(\cA^{(2)}_\proj\bigr)^{sa}$ be of the form
$$
h=Ae+Bf+\psi(efe),
$$
where $\psi$ is a real continuous function on $[0,1]$. Choosing $k(N),l(N)$ such that
$k(N)/N\to\alpha$ and $l(N)/N\to\beta$, we have
\begin{align}\label{F-6.3}
&{1\over N^2}\log\int_{G(N,k(N))\times G(N,l(N))}
\exp\bigl(-N\Tr_N(\psi(P,Q))\bigr)
\,d\bigl(\gamma_{G(N,k(N))}\otimes\gamma_{G(N,l(N))}\bigr)(P,Q) \nonumber\\
&\qquad=-A{k(N)\over N}-B{l(N)\over N}+{1\over N^2}\log\int_{[0,1]^n}
\exp\Biggl(-N\sum_{i=1}^N\psi(x_i)\Biggr)\,d\lambda_N(x_1,\dots,x_N),
\end{align}
where $\lambda_N$ is the empirical eigenvalue distribution of $PQP$ when
$(P,Q)$ is distributed under $\gamma_{G(N,k(N))}\otimes\gamma_{G(N,l(N))}$. By
applying Varadhan's integral lemma (see \cite[4.3.1]{DZ}) to the large deviation in
\cite[Theorem 2.2]{HP1} we have
\begin{align}\label{F-6.4}
&\lim_{N\to\infty}{1\over N^2}\log\int_{[0,1]^n}
\exp\Biggl(-N\sum_{i=1}^N\psi(x_i)\Biggr)\,d\lambda_N(x_1,\dots,x_N) \nonumber\\
&\quad=\sup_\nu\Biggl\{
-(1-\min\{\alpha,\beta\})\psi(0)-\max\{\alpha+\beta-1,0\}\psi(1)
-{1\over2}\int_{[0,1]}\psi(x)\,d\nu(x) \nonumber\\
&\hskip3cm+{1\over4}\Sigma(\nu)
+{|\alpha-\beta|\over2}\int_{[0,1]}\log x\,d\nu(x) \nonumber\\
&\hskip4cm+{|\alpha+\beta-1|\over2}\int_{[0,1]}\log(1-x)\,d\nu(x)-C\Biggr\},
\end{align}
where $\nu$ runs over all measures on $(0,1)$ with $\nu((0,1))=2\rho$. Here, $\rho$ is
in \eqref{F-2.1} and $C$ in \eqref{F-2.2}. For $\tau=\tau_{\nu,\{\alpha_{ij}\}}\in
TS_{(\alpha,\beta)}\bigl(\cA^{(2)}_\proj\bigr)$ (see \S2 and \S5), when
$\alpha_{00}\alpha_{11}=\alpha_{01}\alpha_{10}=0$ (this is necessary for
$\chi_\proj(\tau)>-\infty$), $\chi_\proj(\tau)$ is given as in Proposition \ref{P-2.1}
and moreover we get
\begin{align*}
\tau(h)&=A\alpha+B\beta+(\alpha_{10}+\alpha_{01}+\alpha_{00})\psi(0)
+\alpha_{11}\psi(1)+{1\over2}\int_{(0,1)}(\psi(x)+\psi(0))\,d\nu(x) \\
&=A\alpha+B\beta+(1-\min\{\alpha,\beta\})\psi(0)
+\max\{\alpha+\beta-1,0\}\psi(1)+{1\over2}\int_{(0,1)}\psi(x)\,d\nu(x)
\end{align*}
thanks to \eqref{F-2.3}. Furthermore, Proposition \ref{P-2.1} implies that
$\chi_\proj(\tau)$ is concave and weakly* upper semi-continuous restricted on
$TS_{(\alpha,\beta)}\bigl(\cA^{(2)}_\proj\bigr)$. Hence we obtain (1) and (2) by
\eqref{F-6.3} and \eqref{F-6.4} together with the formulas of $\chi_\proj(\tau)$ and
$\tau(h)$. Moreover, (3) follows from (2) due to the duality for conjugate functions
(or Legendre transforms). Finally, (4) is obvious from (2) and Remark \ref{R-6.3}.
\end{proof}

Now, we introduce a free entropy-like quantity for tracial states on
$\cA^{(n)}_\proj$ (or for $n$-tuples of projections) via the (minus) Legendre
transform of free pressure. For each $\vec\alpha\in[0,1]^n$ and
$\tau\in TS_{\vec\alpha}\bigl(\cA^{(n)}_\proj\bigr)$ define
$$
\eta_\proj(\tau):=\inf\Bigl\{\tau(h)+\pi_{\vec\alpha}(h):
h\in\bigl(\cA^{(n)}_\proj\bigr)^{sa}\Bigr\}.
$$
Since $\pi_{\vec\alpha}$ is a convex continuous function on
$\bigl(\cA^{(n)}_\proj\bigr)^{sa}$ by Proposition \ref{P-6.2}, the above Legendre
transform is reversed so that for every $h\in\bigl(\cA^{(n)}_\proj\bigr)^{sa}$ we have
$$
\pi_{\vec\alpha}(h)=\sup\Bigl\{-\tau(h)+\eta_\proj(\tau):
\tau\in TS_{\vec\alpha}\bigl(\cA^{(n)}_\proj\bigr)\Bigr\}.
$$
For each $h\in\bigl(\cA^{(n)}_\proj\bigr)^{sa}$ there exists a
$\tau\in TS_{\vec\alpha}\bigl(\cA^{(n)}_\proj\bigr)$ such that
$$
\pi_{\vec\alpha}(h)=-\tau(h)+\eta_\proj(\tau).
$$
This equality condition is a kind of variational principle and such $\tau$ may be
called an equilibrium tracial state associated with $h$ (and $\vec\alpha$).

Moreover, for each $n$-tuple $(p_1,\dots,p_n)$ of projections in $(\cM,\tau)$, we have
$\tau_{(p_1,\dots,p_n)}\in TS\bigl(\cA^{(n)}_\proj\bigr)$ defined by
$\tau_{(p_1,\dots,p_n)}(h):=\tau(h(p_1,\dots,p_n))$, where
$h\in\cA^{(n)}_\proj\mapsto h(p_1,\dots,p_n)\in\cM$ is the $*$-homomorphism sending
$e_i$ to $p_i$, $1\le i\le n$. We define
$$
\eta_\proj(p_1,\dots,p_n):=\eta_\proj(\tau_{(p_1,\dots,p_n)}).
$$
It is obvious by definition that the quantity $\eta_\proj(p_1,\dots,p_n)$ has the same
properties as $\chi_\proj(p_1,\dots,p_n)$ given in Proposition \ref{P-1.2}.

\begin{thm}\label{T-6.5}
Let $p_1,q_1,\dots,p_n,q_n,r_1,\dots,r_{n'}$ be projections in $(\cM,\tau)$.
\begin{itemize}
\item[(1)] $\eta_\proj(p_1,\dots,p_n)\ge\chi_\proj(p_1,\dots,p_n)$.
\item[(2)] If $\{p_1,q_1\}$, $\dots$, $\{p_n,q_n\}$, $\{r_1\}$, $\dots$, $\{r_{n'}\}$
are free, then
$$
\eta_\proj(p_1,q_1\dots,p_n,q_n,r_1,\dots,r_{n'})
=\chi_\proj(p_1,q_1\dots,p_n,q_n,r_1,\dots,r_{n'}).
$$
\end{itemize}
\end{thm}

\begin{proof}
The proof of (1) is similar to that of \cite[Theorem 4.5\,(1)]{Hi}. By \eqref{F-6.2}
and Proposition \ref{P-6.4}\,(3), $\eta_\proj=\chi_\proj$ holds when $n=1$ or $n=2$.
Hence (2) is seen from (1) together with the subadditivity of $\eta_\proj$ and the
additivity of $\chi_\proj$ in Theorem \ref{T-4.1}\,(1).
\end{proof}

\begin{remark}\label{R-6.6}{\rm
It is known in a forthcoming paper \cite{HMU} that if
$\chi_\proj(p_1,\dots,p_n)>-\infty$ and
$\sum_{i=1}^n\min\{\tau(p_i),\tau(\1-p_i)\}>1$ (this forces $n\ge3$), then
$\{p_1,\dots,p_n\}''$ is a non-$\Gamma$ II$_1$ factor. Choose two different
$\tau_1,\tau_2\in TS_{\vec\alpha}\bigl(\cA^{(n)}_\proj\bigr)$ such that
$\chi_\proj(\tau_1)$ and $\chi_\proj(\tau_2)$ are finite and
$\sum_{i=1}^n\min\{\alpha_i,1-\alpha_i\}>1$. For $\tau_0:=(\tau_1+\tau_2)/2$ we get
$\eta_\proj(\tau_0)>-\infty$ by the concavity of $\eta_\proj$ on
$TS_{\vec\alpha}\bigl(\cA^{(n)}_\proj\bigr)$. But $\chi_\proj(\tau_0)=-\infty$ due to
the above mentioned fact. Hence, $\eta_\proj$ and $\chi_\proj$ are not equal
in general.
}\end{remark}

Finally, we note that the definition \eqref{F-6.1} is slightly modified in such a way
that the modified free pressure $\pi^{(2)}_{\vec\alpha}(g)$ is defined for self-adjoint
elements of the minimal $C^*$-tensor product
$\cA^{(n)}_\proj\otimes_{\min}\cA^{(n)}_\proj$ and the modified quantity
$\tilde\eta_\proj(p_1,\dots,p_n)$ induced from $\pi^{(2)}_{\vec\alpha}$ via Legendre
transform is always equal to $\chi_\proj(p_1,\dots,p_n)$. We omit the
details concerned with this modification that are essentially same as \cite[\S6]{Hi}.

\end{document}